\documentclass[12pt,a4paper,twoside]{article}
\usepackage{color}
\usepackage{amsthm,amsmath,amssymb}
\newtheorem{teor}{Theorem}[section]
\newtheorem{defin}[teor]{Definition}
\newtheorem{lemm}[teor]{Lemma}
\newtheorem{osse}[teor]{Remark}
\newtheorem{prop}[teor]{Proposition}
\newtheorem{defi}[teor]{Definition}
\newtheorem{coro}[teor]{Corollary}
\newtheorem{prob}[teor]{Problem}

\newcommand{\bele}{\begin{lemm}\begin{sl}}
\newcommand{\enle}{\end{sl}\end{lemm}}
\newcommand{\bedef}{\begin{defi}\begin{sl}}
\newcommand{\eddef}{\end{sl}\end{defi}}
\newcommand{\bete}{\begin{teor}\begin{sl}}
\newcommand{\ente}{\end{sl}\end{teor}}
\newcommand{\beos}{\begin{osse}\begin{rm}}
\newcommand{\eddos}{\end{rm}\end{osse}}
\newcommand{\bepr}{\begin{prop}\begin{sl}}
\newcommand{\empr}{\end{sl}\end{prop}}
\newcommand{\bepro}{\begin{prob}\begin{rm}}
\newcommand{\empro}{\end{rm}\end{prob}}
\newcommand{\bede}{\begin{defin}\begin{sl}}
\newcommand{\edde}{\end{sl}\end{defin}}
\newcommand{\beco}{\begin{coro}\begin{sl}}
\newcommand{\enco}{\end{sl}\end{coro}}
\newcommand{\disp}{\displaystyle}

\newcommand{\thspace}{\hspace{3mm}}

\newcommand{\quext}{\quad\text}
\newcommand{\qquext}{\qquad\text}
\newcommand{\de}{\partial}

\newcommand{\RR}{\mathbb{R}}
\newcommand{\NN}{\mathbb{N}}

\newcommand{\beeq}[1]{\begin{equation}\label{#1}}
\newcommand{\eddeq}{\end{equation}}

\newcommand{\beeqa}[1]{\begin{eqnarray}\label{#1}}
\newcommand{\eddeqa}{\end{eqnarray}}

\newcommand{\beal}[1]{\begin{align}\label{#1}}
\newcommand{\eddal}{\end{align}}

\newcommand{\bespl}[1]{\begin{split}\label{#1}}
\newcommand{\edspl}{\end{split}}

\newcommand{\bega}[1]{\begin{gather}\label{#1}}
\newcommand{\edga}{\end{gather}}

\newcommand{\beeqax}{\begin{eqnarray*}}
\newcommand{\eddeqax}{\end{eqnarray*}}

\def\qed{\ifmmode 
  \else \leavevmode\unskip\penalty9999 \hbox{}\nobreak\hfill
  \fi
  \quad\hbox{\hskip.5em\vrule width.4em height.6em depth.05em\hskip.1em}}
\def\endproofsym{\qed}
\newcommand{\dimbox}{\hbox{\hskip.5em\vrule width.4em height.6em depth.05em\hskip.1em}}
\renewenvironment{proof}[1][Proof]{\trivlist\item[\hskip\labelsep{\hskip0pt
    {\normalfont\scshape#1.}\hskip .321429\parindent}]\ignorespaces}
{\endproofsym\endtrivlist}
\def\endnobox{\def\endproofsym{}\end{proof}\def\endproofsym{\qed}}

\newcommand{\no}{\nonumber}

\newcommand{\beeqao}{\begin{eqnarray}\no}
\newcommand{\bealo}{\begin{align}\no}
\newcommand{\besplo}{\begin{split}\no}
\newcommand{\begao}{\begin{gather}\no}

\newcommand{\nor}[2]{\|#1\|_{#2}}

\newcommand{\cc}{{\mathfrak c}}

\newcommand{\+}{\hspace{1pt}}
\newcommand{\perogni}{\forall\,}
\newcommand{\esiste}{\exists\,}

\newcommand{\ittu}{\int_t^{t+1}}
\newcommand{\io}{\int_\Omega}

\newcommand{\epsi}{\varepsilon}

\newcommand{\X}{\mathcal{X}}

\newcommand{\lhs}{left hand side}
\newcommand{\rhs}{right hand side}

\DeclareMathOperator{\deriv}{d}

\DeclareMathOperator{\dist}{dist}
\DeclareMathOperator{\dom}{dom}

\DeclareMathOperator{\Id}{Id}

\DeclareMathOperator{\sign}{sign}

\DeclareMathOperator{\loc}{loc}
\DeclareMathOperator{\reg}{reg}

\newcommand{\CZV}{C^0([0,T];V)}

\newcommand{\LDH}{L^2(0,T;H)}

\newcommand{\Tuno}{{\cal T}^1}
\newcommand{\Tp}{{\cal T}^p}
\newcommand{\Tpj}{{\cal T}^{p_j}}
\newcommand{\Tq}{{\cal T}^q}

\let\TeXchi\chi
\def\chi{{\setbox0 \hbox{\mathsurround0pt
$\TeXchi$}\hbox{\raise\dp0 \copy0 }}}

\newcommand{\zzn}{_{0,n}}

\newcommand{\Pn}{(P$_n$)}

\newcommand{\calH}{{\mathcal H}}

\newcommand{\calX}{{\mathcal X}}
\newcommand{\calXreg}{{\mathcal X}_{\reg}}
\newcommand{\calXi}{{\mathcal X}_\infty}

\newcommand{\calT}{{\mathcal T}}
\newcommand{\calG}{{\mathcal G}}

\newcommand{\calA}{{\mathcal A}}
\newcommand{\calF}{{\mathcal F}}
\newcommand{\calE}{{\mathcal E}}
\newcommand{\calM}{{\mathcal M}}

\newcommand{\calS}{{\mathcal S}}

\newcommand{\calV}{{\mathcal V}}

\newcommand{\zzu}{_{0,1}}
\newcommand{\zzd}{_{0,2}}

\newcommand{\barr}{\overline{r}}
\newcommand{\debole}{\rightharpoonup}

\newcommand{\barO}{\overline{\Omega}}

\newcommand{\subr}{\underline{r}}

\newcommand{\dit}{\deriv\!t}
\newcommand{\dis}{\deriv\!s}

\newcommand{\dir}{\deriv\!r}

\newcommand{\ddt}{\frac{\deriv\!{}}{\dit}}

\newcommand{\ui}{u_\infty}

\newcommand{\uut}{u_{1,t}}
\newcommand{\udt}{u_{2,t}}

\newcommand{\calBl}{{\cal B}_\ell}

\numberwithin{equation}{section}
\begin{document}

\title{Attractors for the semiflow associated 
with a class of doubly nonlinear parabolic equations}
\author{Giulio Schimperna\\
{\sl Dipartimento di Matematica, Universit\`a di Pavia}\\
{\sl Via Ferrata, 1}\\
{\sl I-27100 Pavia, Italy}\\
{\tt giusch04@unipv.it}\\
\and
Antonio Segatti\\
{\sl Weierstrass Institute for Applied Analysis and Stochastics}\\
{\sl Mohrenstrasse, 39}\\
{\sl D-10117 Berlin, Germany}\\
{\tt segatti@wias-berlin.de}\\
}
\maketitle

\begin{abstract}
 A doubly nonlinear parabolic equation of the form 
 $\alpha(u_t)-\Delta u+W'(u)= f$, complemented with initial
 and either Dirichlet or Neumann homogeneous boundary
 conditions, is addressed. The two nonlinearities are given
 by the maximal monotone function $\alpha$ and by the derivative
 $W'$ of a smooth but
 possibly nonconvex potential $W$; $f$ is a known source.
 After defining a proper notion of solution
 and recalling a related existence result, 
 we show that from any initial datum emanates {\sl at
 least}\/ one solution which gains further regularity
 for $t>0$. Such {\sl regularizing solutions}\/
 constitute a semiflow $\calS$ for which uniqueness is
 satisfied for strictly positive times and we
 can study long time behavior properties. In particular, we 
 can prove existence of both global and 
 exponential attractors and investigate the structure 
 of $\omega$-limits of single trajectories.
\end{abstract}

\noindent
{\bf Key words:}\thspace doubly nonlinear equation, singular
potential, semiflow, global attractor, energy method,
$\omega$-limit.
  
\vspace{2mm}

\noindent
{\bf AMS (MOS) subject clas\-si\-fi\-ca\-tion:}\thspace
35K55, 35B40, 35B41.

\vspace{2mm}



\pagestyle{myheadings}
\newcommand\testopari{\sc G.~Schimperna -- A.~Segatti}
\newcommand\testodispari{\sc Doubly Nonlinear Equations}
\markboth{\testodispari}{\testopari}


\section{Introduction}
In this paper we are interested in the following doubly non linear
 parabolic equation
\begin{equation}
  \alpha(u_t)-\Delta u+W'(u)=f,\quad\mbox{for a.e.~}\,(x,t)\in 
   \Omega\times (0,+\infty),\label{eq-concreta}
\end{equation}
where $\Omega\subset \mathbb{R}^N$,
$1\le N\le 3$, is a bounded domain with smooth
boundary $\partial\Omega$. Here $\alpha$ is a differentiable and strongly 
monotone (i.e., $\alpha'\ge\sigma>0$) function in $\mathbb{R}$, 
$W'$ is the derivative of a $\lambda$-convex 
(i.e., $W''\ge-\lambda$, $\lambda\ge 0$)
configuration potential, and $f$ is a source.
The equation is complemented with the initial conditions
and with homogeneous boundary conditions of either Dirichlet or
Neumann type. Equations like \eqref{eq-concreta},
apart from their own mathematical interest, 
can arise in large variety of applications, as the modelization of 
phase change phenomena \cite{BDG,BFL,FV,LSS,LSS2}, gas flow
through porous media \cite{Fr}
and damaging of materials \cite{BS,FKS,MR}.

Existence of (at least) one solution to 
initial-boundary value problems
for a class of doubly nonlinear equations
including~\eqref{eq-concreta} was proved
in the paper~\cite{CV} 
(see also \cite{Arai,Barbu75,Senba} 
for preceding related results).
The questions of regularity, uniqueness,
continuous dependence on data
and long time behavior of solutions, however, 
were not dealt with in~\cite{CV} and remained widely open
for a long time. Moreover, the results of
\cite{CV} require the restrictive assumption that
$\alpha$ is bounded in the sense of operators
(i.e.~it maps bounded sets into bounded sets),
which is not always fulfilled in physical applications
(see the papers quoted above referring to 
specific models). On account of these considerations,
in the former paper \cite{SSS1}, written in collaboration
with U.~Stefanelli, we 
introduced a 
new concept of solution ({\sl stronger}\/ than
that in \cite{CV}, see Def.~\ref{defisolX2} below) 
and showed existence of this kind of solution with essentially 
no restriction on $\alpha$. This permitted to prove also 
uniqueness, at least in some special cases, as well as
existence of nonempty $\omega$-limits.
A further contribution in this field
has been recently given in
\cite{EZ}, where a doubly nonlinear equation
strictly related to~\eqref{eq-concreta}, but of
{\sl degenerate}\/ type, is addressed from 
the viewpoint of both well-posedness and long time behavior.

One of the main issues of this paper is
a regularization property, 
holding for $t>0$, of the solutions to the IBV problem for 
\eqref{eq-concreta}. 
Due to the strong parabolicity of the system 
($\alpha'\ge\sigma>0$) such a fact is to 
be expected; however, the proof requires a somehow
tricky machinery due to the presence of very general 
nonlinearities. The key point, resembling in some way
the approach given also in \cite{EZ}, consists in 
an Alikakos-Moser \cite{Al} iteration scheme,
operated here on the (formal) time derivative 
of~\eqref{eq-concreta}, coupled with the 
use (infinitely many times) of the {\sl uniform}\/
Gronwall lemma (see, e.g., \cite{Te}).
In this way we demonstrate that, 
if the source $f$ is essentially bounded, then
there exist solutions $u(t)$ (called ``regularizing
solutions'' in the sequel, see 
Def.~\ref{defisolXi}) which, for $t>0$,
are in $L^\infty(\Omega)$ together with their Laplacian 
and with $W'(u(t))$. Moreover, for $t>0$ uniqueness 
holds, whereas from any initial
datum can start more than one trajectory
unless the datum is more regular itself. 

The regularization property serves also as a starting
point to improve the results of \cite{SSS1} 
regarding long time behavior. 
Actually, in case the potential $W$ is {\sl real
analytic}\/ we can show, using the {\sl Simon-\L ojasiewicz}\/
method (cf.~\cite{Lo1,Lo2,Si2}, see also 
\cite{CJ,GPS2}), that $\omega$-limits of all
single trajectories contain only one point. This can be done
without the severely restrictive assumptions on 
the growth of $\alpha$ at $\infty$ which were considered
in \cite{SSS1}. We remark that the Simon-\L ojasiewicz
method is a deep and powerful 
tool that in recent year has been
applied to characterize $\omega$-limit 
sets of solutions to several different types
of nonlinear evolution equations (see, e.g., 
\cite{Ch,CJ,FS,HJ,Je} among the 
many related works).

From the viewpoint of long time behavior,
however, our main result regards the existence and regularity 
properties of attractors. We notice that
a contribution to this analysis has been recently 
given in~\cite{Se}, where a (rather weak) notion of 
global attractor is introduced for
a class of equations including \eqref{eq-concreta}.
However, due to the very general and abstract setting 
adopted there (very similar to that of \cite{CV}),
the attractor constructed in \cite{Se} seems
not very flexible from the point 
of view of regularity (more precisely, it appears
difficult to characterize it beyond
its mere existence property). Moreover, 
the result in \cite{Se} holds only
under the boundedness assumption on $\alpha$ considered
in \cite{CV} and consequently is not suitable
for our specific situation.

Here, also thanks to the much more specific form
of equation \eqref{eq-concreta},
we can prove the existence of a global
attractor in the natural phase space (i.e.~under the
precise conditions ensuring existence). The key
point is the use of the so-called 
{\sl energy method}\/ by J.~Ball 
(cf.~\cite{Ba2}, see also \cite{moise-rosa-wang04} 
and the references therein),
which permits to prove this result without
reinforcing the conditions on the source $f$ (namely,
we do not need to ask summability of its space
derivatives) and despite the apparent lack
of a dissipative estimate in the natural phase
space (see Remark~\ref{ossedissi} below). 
We point out that,
due to the (possible) non-uniqueness at $t=0$,
the semiflow $\calS$ associated to 
\eqref{eq-concreta} for which we can prove
existence of the global attractor has to be
carefully defined (in particular, ``nonregularizing'' 
solutions have to be excluded, see Remark~\ref{ancorarego}).
This is in agreement with other works where 
equations with (at least partial) lack of uniqueness
are addressed (see, e.g., \cite{Ba1,Ba2,MP,RSS,Se,Se2}).

Our final issue is concerned with exponential attractors,
whose existence is proved by using as a technical tool
the so-called method of {\sl short trajectories}\/
(or {\sl $\ell$-trajectories}) due to M\'alek and 
Pra\v z\'ak \cite{MP}. Actually, this device
permits to get in a simple way the contractive estimates
required to have the exponential attraction property.
We stress that this approach is quite similar to that used in
\cite{Mi}, where the equation (strictly related to
\eqref{eq-concreta} or, more precisely, to its time
derivative)
\begin{equation}
  \disp \alpha(u)_t-\Delta u+W'(u)=f,\quad\mbox{for a.e.~}\,(x,t)\in 
   \Omega\times (0,+\infty),\label{eq-miranv}
\end{equation}
is addressed (although under partly different assumptions on
the nonlinearities).

We conclude with the plan of the paper. In the next Section
some preliminary material is recalled. Next, our results
are presented in a rigorous way in Section~3, where in particular
the required notions of solution are introduced. The subsequent
Section~4 contains the proof of the regularization 
property and Sections~5 and~6 are devoted to global
and exponential attractors, respectively. Finally,
an abstract existence Theorem for global attractors,
partially generalizing \cite[Thm.~3.1]{Ba1}, is reported
in the Appendix.


\section{Preliminaries}\label{preli}

In this section we introduce some notations and recall some
preliminary notions which are needed 
to state our problem in a precise way.
First of all, we set $H:=L^{2}(\Omega)$ and denote by
$(\cdot,\cdot)$ the scalar product in $H$ and by
$\|\cdot\|$ the related norm. The symbol $\|\cdot\|_{X}$ 
will indicate the norm in the generic Banach space $X$. 
Moreover, focusing on the Dirichlet case, we set
$V:=H^{1}_{0}(\Omega)$, $V':=H^{-1}(\Omega)$ and 
identify $H$ and $H'$ so that we obtain the Hilbert 
triplet $V\subset H\subset V'$, where inclusions
are continuous and compact.
The notation $\langle\cdot,\cdot\rangle$
will stand for the duality between $V'$ and $V$. 
We also let $B:V\to V'$ denote the distributional Laplace operator,
namely
\beeq{defiB}
  B:V\to V', \quad\langle Bu,v\rangle=(\nabla u,\nabla v)~~
   \perogni u, v\in V.
\end{equation}
\beos\label{dirineum}
 Here and in the sequel we assumed Dirichlet conditions just
 for simplicity. Indeed, the (homogeneous) Neumann case works 
 as well with the following simple change:
 we have to set $V:=H^1(\Omega)$, $V':=H^1(\Omega)'$
 and, in place of \eqref{defiB},
 \beeq{defiBneum}
  B:V\to V', \quad\langle Bu,v\rangle=(u,v)+(\nabla u,\nabla v)
   ~~\perogni u, v\in V.
 \end{equation}
 All the results and proofs in the sequel then still 
 work with no further change.
\eddos
\noindent%
In order to correctly describe the asymptotic behavior of solutions,
we need to introduce the space of
{\sl $L^p_{\loc}$-translation bounded}\/ functions.
As $X$ is a Banach space and $p\in[1,+\infty)$ we set
\beeq{defiTp}
  \Tp(T,\infty;X):=\Big\{v\in L^p_{\loc}(T,\infty;X):
   \sup_{t\ge T}\ittu\nor{v(s)}{X}^p\,\dis<\infty\Big\},
\end{equation}
which is a Banach space with respect
to the natural (graph) norm
\beeq{defiTpnorm}
  \nor{v}{\Tp(T,\infty;X)}^p
   :=\sup_{t\ge T}\ittu\nor{v(s)}{X}^p.
\end{equation}
Next, we recall the {\sl uniform}\/ Gronwall 
Lemma (see, e.g., \cite[Lemma~III.1.1]{Te}),
which will be repeatedly used in the sequel:
\bele\label{unigronw}
 Let $y,a,b\in L^1_{\loc}(0,+\infty)$ three non negative
 functions such that $y'\in L^1_{\loc}(0,+\infty)$ and,
 for some $T\ge 0$, 
 \beeq{disunigronw} 
   y'(t)\le a(t)y(t)+b(t)
    \qquext{for a.e.~\,$t\ge T$,}
 \end{equation}
 and let $k_1,k_2,k_3$ three nonnegative constants such
 that 
 \beeq{norTpgronw} 
   \nor{a}{\Tuno(T,\infty;\RR)}\le k_1, \qquad
   \nor{b}{\Tuno(T,\infty;\RR)}\le k_2, \qquad
   \nor{y}{\Tuno(T,\infty;\RR)}\le k_3.
 \end{equation}
 Then, we have that 
 \beeq{tesiunigronw} 
   y(t+\tau)\le \big(k_2+k_3/\tau\big)e^{k_1}
    \qquext{for all~\,$t\ge T$}.\dimbox
 \end{equation}
\enle
\noindent%
Now, let us recall some basic facts about absorbing sets and
attractors. Assuming that $\calX$ is a complete metric 
space, we shall (conventionally) call
a {\sl semiflow}\/ on $\calX$ a family
$\calS$ of maps from $[0,\infty)$ to $\calX$,
called {\sl trajectories}, or {\sl solutions},
satisfying properties 
{\sl (S1)-(S5)}\/ listed below. We stress that
this definition, which partly follows the approach in
\cite{Ba1,Ba2} (see also \cite{RSS}),
is not standard at all. Actually,
in Ball's terminology, $\calS$ could be noted like
a ``strongly-weakly continuous generalized semiflow with
unique continuation''. We say here
``semiflow'' just for brevity.\\[2mm]
{\sl (S1 -- existence)}~~For all $u_0\in \calX$ there exists 
 at least one $u\in \calS$ such that $u(0)=u_0$;\\[1mm]
{\sl (S2 -- translation invariance)}~~For all $u\in \calS$
 and $T\ge 0$, the map $v:[0,\infty)\to \calX$ given by
 $v(t):=u(T+t)$ still belongs to $\calS$;\\[1mm]
{\sl (S3 -- concatenation)}~~For all $u,v\in \calS$ such 
 that for some $T>0$ it is $u(T)=v(0)$, the map 
 $z:[0,\infty)\to \calX$ coinciding con $u$ on
 $[0,T]$ and given by $z(t)=v(t-T)$ on $(T,\infty)$ 
 belongs to $\calS$;\\[1mm]
{\sl (S4 -- unique continuation for $T>0$)}~~
 For all $u,v\in \calS$ such that $u(T)=v(T)$ for some 
 $T>0$, it is $u(t)=v(t)$ for all $t\in[T,\infty)$;\\[1mm]
{\sl (S5 -- strong-weak semicontinuity)}~~
 We assume that, beyond the {\sl strong}\/ topology 
 induced by the metric, $\calX$ is endowed with a {\sl weaker}\/
 topology. Then, we firstly 
 ask that all elements of $\calS$ are
 {\sl weakly continuous}\/ from $[0,\infty)$ to
 $\calX$. Next, that for all sequence $\{u_n\}\subset\calS$ such
 that $u_n(0)=:u\zzn$ tends {\sl strongly}\/ (i.e.~with 
 respect to the metric) to some $u_0\in\calX$, there 
 exist a subsequence (not relabelled) of $\{u_n\}$ 
 and $u\in\calS$ with $u(0)=u_0$ such that,
 for all $t>0$, $u_n(t)$ tends {\sl weakly}\/ to $u(t)$.
\beos\label{sullacont}
 Regarding {\sl (S5)}, if $\calX$ is a Banach space, a
 natural choice for the ``weak topology'' mentioned
 there is of course that induced by the weak 
 (or, in some cases, the weak star) convergence. 
 We will show in the sequel (see in particular the 
 Appendix) that the lack of a more usual 
 ``strong-strong'' continuity property does not
 prevent use of time regularization-compactness
 methods to get existence of the global 
 attractor. This fact has been noted also in 
 other recent papers \cite{PZ,ZYS}.
\eddos
\noindent%
We assumed property {\sl (S4)}, which is not completely
standard, just to fit the case of our system for which
uniqueness holds only from $t>0$.
If $\calS$ is a semiflow, we define the {\sl space of 
regularized values} of $\calS$ as
\beeq{defiV}
  \calXreg:=\big\{u(t): u\in\calS,t>0\big\}.
\end{equation}
Moreover, if $u\in \calS$, we recall that the (strong)
$\omega$-limit of $u$ is the set of all limit 
(w.r.t.~the metric) 
points of subsequences of $u(t)$ as $t\nearrow\infty$.
From {\sl (S2)}\/ and {\sl (S4)}, it is apparent that
it can be naturally associated to a semiflow $\calS$
the family $\{S(t)\}$, $t\in[0,\infty)$,
of operators from $\calXreg$ to itself,
with $S(t)$ mapping $x\in\calXreg$ into $u(t)$, 
where $u\in\calS$ is the (unique) trajectory
such that $u(0)=x$. It is then clear
that $\{S(t)\}$ satisfies the usual semigroup 
properties. Due to the lack of uniqueness, $S(t)$ cannot
be extended to the whole $\calX$. Nevertheless, we can 
introduce the family of {\sl multivalued}\/
mappings $\{T(t)\}$, $t\in[0,\infty)$, given by
\beeq{defiTt}
  T(t):\calX\to 2^{\calX},\qquad
   T(t)u:=\big\{v(t):v\in\calS,v(0)=u\big\}
\end{equation}
and by {\sl (S4)} it is then clear that the restriction of 
$T(t)$ to $\calXreg$ coincides with $S(t)$.

Next, we recall that a compact subset 
$\mathcal{A}$ of the {\sl phase space}\/
$\mathcal{X}$ is the {\sl global attractor}\/ for the semiflow 
$\calS$ if the following conditions are satisfied:\\[2mm]
{\sl (A1)}~~The set $\mathcal{A}$ is fully invariant, i.e.,
 $T(t)\mathcal{A}=\mathcal{A}$ for all $t\ge 0$;\\[1mm]
{\sl (A2)}~~$\mathcal{A}$ attracts the images of all
 bounded subsets of $\mathcal{X}$ as $t\nearrow +\infty$, 
 namely
 \begin{equation}\label{attrazione}
   \lim_{t\nearrow+\infty}\dist(T(t)B,\mathcal{A})=0,
    \mbox{~~for all bounded }\,B\subset \mathcal{X},
 \end{equation}
 where $\dist$ is the standard {\sl non-symmetric}\/
 Hausdorff distance between sets in
 $\mathcal{X}$.\\[2mm]
We point out that the global attractor
represents the first (although extremely important)
step in understanding the long-time dynamics of a given
evolutive system.
However, it may also present some drawbacks. 
First of all, it may be reduced 
to a single point, thus failing in capturing all the transient behaviour of
the process. Moreover, in general it is extremely difficult 
to estimate the rate of convergence in \eqref{attrazione} and to 
express it in terms of physical parameters. 
In this regard, simple examples 
show that this rate of convergence 
may be arbitrarily slow. This fact makes the global attractor
very sensitive to perturbations and to numerical approximation.
The concept of {\sl exponential attractor}\/ has then been proposed
(see, e.g., \cite{EFNT}) to possibly overcome 
this difficulty. We recall that a compact
subset $\mathcal{M}$ of the phase space $\calX$ is called 
an {\sl exponential attractor}\/ for the semiflow $\calS$ 
if the following conditions are satisfied:\\[2mm]
{\sl (E1)}~~The set $\mathcal{M}$ is positively invariant, i.e.,
 $T(t)\mathcal{M}\subset \mathcal{M}$ for all $t\ge 0$;\\[1mm]
{\sl (E2)}~~The fractal dimension (see, e.g., \cite{Ma})
 of $\mathcal{M}$ in $\mathcal{X}$ is finite;\\[1mm]
{\sl (E3)}~~The set $\mathcal{M}$ attracts exponentially 
 fast the images of the bounded sets $B$ of the phase space
 $\mathcal{X}$. Namely, for
 every bounded $B\subset \mathcal{X}$ there exist $C,\beta>0$ 
 depending on $B$ and such that 
 \begin{equation}\label{expo-attractio}
   \dist(T(t)B,\mathcal{M})\le C\+e^{-\beta t},
   \quad\perogni t\ge 0.\\[1mm]
 \end{equation}
Thanks to $(E3)$ it follows that, compared to the global 
attractor, an exponential attractor
is much more robust to perturbation and to the important 
issue of numerical approximation (see, e.g., \cite{EFNT} and
\cite{FaGa}).
Moreover, when the exponential attractor $\mathcal{M}$ exists, it
contains the global attractor $\mathcal{A}$. Thus, in this case 
also $\calA$ has finite fractal dimension. We point out that, however,
also the theory of exponential attractors presents some disadvantages,
like the lack of uniqueness of $\calM$, whose choice or
construction may be in some sense artificial. 
However, we refer to \cite{Efen-Mir-Zelik06} 
where it is proposed a construction of an exponential attractor 
which selects a proper one valued branch of the exponential 
attractors depending in an H\"older continuous way
on the dynamical system under study.
In recent years several different techniques have been provided to
guarantee existence of exponential attractors. Beyond the original
method \cite{EFNT} based on a direct verification of 
the {\sl discrete squeezing property},
we quote the ``decomposition technique'' developed in \cite{EMZ}
and, in particular, the so-called method of ``$\ell$-trajectories''
(or ``short trajectories''), introduced by M\'alek and Pra\v z\'ak
in \cite{MP}, which provides a simplified framework 
which can be adopted to verify the theoretical conditions 
of \cite {EFNT} leading
to existence of ${\calM}$. Since we shall use this method in the
sequel, we recall here, for convenience of the reader,
its highlights, partly adapting the presentation
in \cite{MP} to our more specific framework.

Let $\X$ be a {\sl Hilbert}\/ space
%
%
%
and, for given $\tau>0$, let us set 
$\X_\tau:= L^2(0,\tau,\X)$.
We assume that there exists a subset 
$B_1$ of $\calX$ such that for any $u_0\in B_1$ there
exists at least one map $u\in C_w([0,\infty);\calX)$
such that $u(0)=u_0$. These maps $u$ are called ``solutions''
in what follows, and we assume that they form a semiflow
$\calS$ on the set $B_1$ endowed with the strong 
and weak topologies inherited from $\calX$. 
We then introduce the space of $\ell$-trajectories 
(where $\ell>0$) as
\beeq{defiltraj}
  \calBl^1:=\left\{\chi:(0,\ell)\to \X, \;\;\chi \mbox{ is a
     solution on } (0,\ell)\right\}.
\end{equation}
The space $\calBl^1$ inherits its topology from $\X_\ell$.
Moreover, according to {\sl (S4)}, any $\ell$-trajectory has, 
among all solutions, unique continuation. We shall assume that 
\beeq{calBl1}
  \calBl^1~~\text{is relatively compact in }\,\X_\ell.
\end{equation}
Then, the method of $\ell$-trajectories basically consists 
in lifting the dynamical system from the phase 
space of initial conditions to the space $\calBl^1$
of $\ell$-trajectories. In particular, by {\sl (S4)} we can 
define a semigroup $L_t$ on $\calBl^1$ by setting
\begin{equation}\label{semigruppoL}
  \left\{L_t\chi\right\}(\tau):=u(t+\tau), \quad \tau\in [0,\ell],
\end{equation}
where $\chi$ is an $\ell$-trajectory and 
$u$ is the {\sl unique}\/ solution such that $u|_{[0,\ell]}=\chi$. 
Then, the assumptions that lead to the existence of the exponential
attractor in the space of $\ell$-trajectories endowed
with the topology of $\calX_\ell$ read as follows 
(see \cite{MP}):\\[2mm]
 {\sl (M1)}~~There exist a space $W_\ell$ compactly embedded 
 into $\X_\ell$ and $\tau>0$ such that
 $L_\tau:\X_\ell\to W_\ell $ is Lipschitz 
 continuous on $\calBl^1$;\\[1mm]
 {\sl (M2)}~~For all $\tau>0$ the family of operators
  $L_t:\X_\ell\to \X_\ell$ is uniformly
 (w.r.t.~$t\in [0,\tau]$) Lipschitz continuous on $\calBl^1$;\\[1mm]
 {\sl (M3)}~~For all $\tau>0$ there exist $c>0$ and
 $\beta \in (0,1]$ such that for all
 $\chi\in \calBl^1$ and $t_1,t_2\in [0,\tau]$ it holds 
 that $\| L_{t_1}\chi-L_{t_2}\chi\|_{\X_\ell}
    \le c\vert t_1-t_2\vert^{\beta}$.\\[2mm]
In \cite[Theorem 2.5]{MP} it is proved that, under the
assumptions above, there exists an exponential attractor 
$\mathcal{M}_\ell$ for the dynamical system $L_t$ on $\calBl^1$.
One of the striking features of this method is that, 
once we have constructed an exponential attractor
in the space of $\ell$-trajectories, we can recover 
the dynamics in the original phase space $B_1$ and obtain
an exponential attractor $\mathcal{M}$ for the semiflow $\calS$. 
To this end, we introduce the {\sl evaluation
map}\/ assigning to a given $\ell$-trajectory 
$\chi$ its end point, i.e.,
\begin{equation}\label{evaluation}
  e:\calBl^1\to \X,\quext{given by }\, e(\chi):=\chi(\ell).
\end{equation}
If it additionally holds that\\[2mm]
{\sl (M4)}~~The map $e$ is H\"older continuous on $\calBl^1$,\\[2mm]
we then obtain the exponential attractor in the phase space
as the image of $\mathcal{E}_\ell$ (see \cite[Theorem 2.6]{MP}), 
namely we have that $\mathcal{M}:=e(\mathcal{M}_\ell)$ is an exponential 
attractor for the semiflow $\calS$ on the space $B_1$.
\beos\label{nuovaespo}
 In general, the semiflow $\calS$  is originally defined 
 on a space ``larger'' than the bounded
 set $B_1$ (usually, but not in our case, on the whole $\calX$),
 and $B_1$ is chosen ``a posteriori'' as a bounded, absorbing
 and positively invariant set
 for the ``original'' $\calS$. One of the advantages 
 of this approach is then that property
 \eqref{calBl1} requires in general very little smoothing
 effect (and is usually straightforward to be checked
 in concrete situations). We also note that, once 
 we have the exponential attractor $\calM$ on $B_1$ 
 where $B_1$ is absorbing, then $\calM$ turns out to be 
 an exponential attractor on the whole space.
\eddos


\section{Main results}

We begin by specifying our basic assumptions on data. First of all,
we ask 
\begin{equation} \label{alfadalbasso}\tag{hp$\alpha$}
   \alpha\in C^1(\RR;\RR),~~\alpha(0)=0,~~
   \alpha'(r)\ge \sigma>0  \mbox{~~for all }\, r\in \mathbb{R}.
\end{equation}
Next, given $\lambda\ge 0$ and
an {\sl open}\/ (either bounded or unbounded) interval 
$I\subset \RR$ with 
$0\in I$, we assume that the potential $W$ fulfills

\beal{W}\tag{hpW1}
  & W\in C^{1,1}_{\loc}(I;\RR),\quad W'(0)=0, \quad
   W''\ge -\lambda\text{~~a.e.~in~}\,I,\\
 \label{coercW}\tag{hpW2}
  & \lim_{r\to\de I}W'(r)\sign r=+\infty.
\end{align}
Property \eqref{W} is called {\sl $\lambda$-convexity}\/
in what follows (see \cite{AGS} for the definition).
Since $W$ is defined up to an additive constant, it 
is also not restrictive to suppose that
\begin{equation} \label{W3}
  \esiste \eta>0:~~W(r)\ge \frac{\eta r^2}{2}
   \mbox{~~for all }\, r\in I.
\end{equation}
We then introduce the basic {\sl phase space}\/ for our
analysis:
\begin{equation} \label{defiX2}
  \calX_2:=\big\{u\in H:Bu,\,W'(u)\in H\big\},
\end{equation}
which is endowed with the metric
\begin{equation} \label{defidX2}
  d_2^2(u,v):=\|u-v\|^2+\|Bu-Bv\|^2
   +\|(W'+\lambda)(u)-(W'+\lambda)(v)\|^2.
\end{equation}
Proceeding as in cite \cite[Lemma~3.8]{RS}
(compare also with \cite[Sec.~3]{Se}),
it is easy to show that $\calX_2$ is a {\sl complete}\/
metric space. It is also clear that $\calX_2\subset V\cap H^2(\Omega)$
(continuously); however, if $I\not=\RR$, in general the inclusion is strict.

We can now list our hypotheses on the initial and source data:
\beal{regodata}\tag{hp$u_0$}
  & u_0\in\mathcal{X}_2,\\
 \label{regof}\tag{hp$f$}
  & f\in L^{\infty}(\Omega).
\end{align} 
%
Then, standardly identifying $\alpha$ and $W'$ as operators from $H$ 
to itself, we introduce the 
\bede\label{defisolX2}
 We call an {\rm $\calX_2$-solution} to the Problem~{\rm (P)}
 given by
 \beal{eqn}
   & \alpha(u_t)+B u+W'(u)= f, \quext{in }\,H,
    \quext{a.e.~in~}\,(0,\infty),\\
  \label{iniz}
   & u|_{t=0}=u_0, \quext{in }\,H
 \end{align}
 one function $u:[0,\infty)\to H$ satisfying\/ 
 \eqref{eqn}, \eqref{iniz} and, for some $C>0$,
 \beal{regou} 
   & u,\, u_t,\, \alpha(u_t),\, Bu,\, W'(u) \in L^\infty(0,\infty;H),\\
  \label{regou2}
   & d_2^2(u(t),0)=\|u(t)\|^2+\|Bu(t)\|^2+\|(W'+\lambda)(u(t))\|^2
       \le C^2\quext{for {\sl all}~}\,t\in[0,\infty).
 \end{align}
\edde
\noindent%
We note that \eqref{eqn}--\eqref{iniz} give a rigorous formulation
of the IBV problem for \eqref{eq-concreta}. With 
condition~\eqref{regou2} we ask the solution to stay 
in the phase space $\calX_2$ for any (and not just a.e.) 
value of the time variable. We can now recall the statement
of the existence result proved in \cite[Thm.~2.5]{SSS1}:
\bete\label{teoexist}
 Assume\/ \eqref{alfadalbasso}, \eqref{W}--\eqref{coercW}, and
 \eqref{regodata}--\eqref{regof}. More precisely,
 suppose that for some $\kappa>0$ it is 
 \beeq{inizk}
   d_2^2(u_0,0)=\|u_0\|^2+\|Bu_0\|^2+\|(W'+\lambda)(u_0)\|^2
       \le \kappa^2.
 \end{equation}
 Then, Problem~{\rm (P)} admits at least
 one $\calX_2$-solution, which additionally satisfies
 \beeq{stimaL2H}
   \|u_t\|_{L^2(0,t;V)}^2\le C^2.
 \end{equation}
 Moreover, the constants $C$ in\/ \eqref{regou2},
 \eqref{stimaL2H} depend only on $\Omega$, $\alpha$, $W$, $f$, 
 and (linearly) on $\kappa$ in\/ \eqref{inizk}. In particular, 
 they do not depend on the time~$t$.
\ente
\noindent%
We remark that~\eqref{stimaL2H}, which was not stated
in \cite[Thm.~2.5]{SSS1} since the coercivity hypotheses
on $\alpha$ considered there were weaker, follows easily
from the proof in \cite[Sec.~3]{SSS1} thanks to the
last assumption in~\eqref{alfadalbasso}.
Let us now see that {\sl some}\/
solutions to Problem~(P)
gain more spatial regularity for $t>0$. With this aim,
we introduce the new space 
\begin{equation} \label{defiXi}
  \calX_\infty:=\big\{u\in L^\infty(\Omega):Bu,\,W'(u)\in L^\infty(\Omega)\big\},
\end{equation}
which is naturally endowed with the (complete) metric
\begin{equation} \label{defidXi}
  d_\infty^2(u,v):=\|u-v\|^2_{L^\infty(\Omega)}+\|Bu-Bv\|_{L^\infty(\Omega)}^2
   +\|(W'+\lambda)(u)-(W'+\lambda)(v)\|_{L^\infty(\Omega)}^2.
\end{equation}
We also introduce weaker notions
of convergence (and, in fact, weaker topologies)
on the spaces $\calX_2$, $\calXi$.
Namely, we say that a sequence $\{u_n\}$ tends
to $u$ weakly in $\calX_2$ (in $\calXi$)
if $u_n\to u$, $Bu_n\to Bu$, and
$(W'+\lambda)(u_n)\to (W'+\lambda)(u)$
weakly in $H$ (weakly {\sl star}\/ in $L^\infty(\Omega)$,
respectively). When we construct below the semiflow $\calS$
on $\calX_2$, property {\sl (S5)} will be implicitly intended
with respect to this weak structure.

To proceed, we need to introduce a 
couple of functionals defined on the space $\calX_2$,
the first of which has the meaning of {\sl energy}:
\beal{defiE}
  \calE(u) & :=\io\Big[\frac{|\nabla u|^2}2
   +W(u)-fu\Big],\\
 \label{defiF}
  \calF(u) & :=\frac12\|Bu+W'(u)\|^2
   -(f,Bu+W'(u)).
\end{align}
It is clear that, since \eqref{W3} and \eqref{regof} hold,
both functionals are finite and bounded from below on 
$\calX_2$. Moreover, mimicking the procedure given in
\cite[Sec.~3]{SSS1}, i.e., {\sl formally}\/ testing
\eqref{eqn} by $\lambda u_t + (Bu+W'(u))_t$, 
and using in particular \eqref{W}, one can expect
that solutions $u$ to Problem~(P) satisfy
\beeq{formLiap}
  \ddt\big(\lambda\calE+\calF\big)(u(t))
   \le 0\quext{for a.e.~}\/t\ge 0.
\end{equation}
Setting then $\calG:=\lambda\calE+\calF$
and noting that there exist $\eta_1,\eta_3>0$ and $\eta_2\ge 0$ such
that
\beeq{Liapd2}
  \eta_1 d_2^2(u,0)-\eta_2
    \le \calG(u)
    \le \eta_3 \big(d_2^2(u,0)+1\big) \quad\perogni u\in \calX_2,
\end{equation}
relation \eqref{formLiap} takes the form of a decay
(or Liapounov) condition for the distance $d_2$.

However, the formal procedure used to get 
\eqref{formLiap} seems very difficult to be justified
if we just know that $u$ is an $\calX_2$-solution.
Actually, \eqref{eqn} is settled in $H$ and
\eqref{regou} does not imply that the test function
$\lambda u_t + (Bu+W'(u))_t$ takes values in $H$.

To overcome this difficulty, we recall that the existence 
Theorem~\ref{teoexist} was shown in \cite{SSS1}
via approximation and compactness methods. We sketch
here, and partly refine, just the highlights
of this procedure. Let us
substitute $\alpha$ and $W$ in \eqref{eqn}
with regularized functions $\alpha_n$ and $W_n$
still satisfying \eqref{alfadalbasso},
\eqref{W} and such that 
\beal{regon1}
  & \alpha_n,(W_n'+2\lambda\Id)\,
   \text{ are Lipschitz continuous with their inverses},\\
 \label{regon2}
  & \alpha_n,(W_n'+2\lambda\Id)\to\alpha,(W'+2\lambda\Id)\,
   \text{ in the sense of graphs \cite{At}},
\end{align}
the latter convergences intended as $n\nearrow\infty$.
Then, noting as \Pn\ the problem still
given by \eqref{eqn} (with the regularized functions)
and \eqref{iniz} (note that the initial datum
is {\sl not}\/ regularized), it is not difficult to
show the
\bepr\label{daPnaP}
 For every $n>0$, 
 Problem~{\rm \Pn}\ has one and only one solution
 $u_n$ such that 
 \beeq{regoun}
   u_{n,tt}\in L^2(0,\infty,H), \qquad 
    u_n,u_{n,t}\in L^2(0,\infty,H^2(\Omega)).
 \end{equation}
 Moreover, $u_n$ satisfies estimates 
 \eqref{regou}, \eqref{regou2}
 with $C$ independent of $n$. Finally,
 for any subsequence of $n\nearrow\infty$, 
 there exists a subsubsequence (still noted here
 as $u_n$) such that $u_n$ suitably (i.e., in the sense 
 specified by \eqref{regou} and \eqref{regou2})
 tends to $u$, where $u$ is an $\calX_2$-solution
 to Problem~{\rm (P)}.
\empr
We point out that the proof of the above Proposition
could be performed just by refining the estimates 
and the passage to the limit in \cite[Sec.~3]{SSS1}.
We omit, for brevity, the technical details of the argument
and rather focus our attentions on the more subtle
consequences of working with solutions 
$u_n$ of~\Pn. Clearly, the functions $u_n$ 
do satisfy \eqref{formLiap} (where, of course, $W_n$ 
replaces $W$ in $\calG$). However, the convergence
$u_n\to u$ specified by estimate \eqref{regou2}
is too weak to let \eqref{formLiap} pass to 
the limit with $n$. Moreover, due to nonuniqueness
for the problem (P), there might exist some
$\calX_2$-solutions which are not, or at least are 
not known to be, limit of (sub)sequences
of solutions to (P$_n$). Actually, we
shall note in the sequel
as {\sl limiting}\/ (respectively, 
{\sl nonlimiting}\/) the solutions to (P)
which are (respectively, are not) limits of 
(sub)sequences of solutions to~\Pn.
%
%
We then introduce a new concept of solution, where
a (much weaker than \eqref{formLiap}) 
form of Liapounov property (cf.~\eqref{intLiap} below)
for $\calG$ is {\sl postulated}.
From the proofs, it will be clear that 
all limiting solutions fulfill \eqref{intLiap},
but there might also exist nonlimiting solutions
satisfying it.
\bede\label{defisolXi}
 A\/ {\rm regularizing solution} to 
 Problem~{\rm (P)} is an
 $\calX_2$-solution which, additionally, fulfills
 the\/ {\rm regularization property}
 \beeq{regoui} 
    u_t,\, \alpha(u_t),\, Bu,\, W'(u) \in L^\infty(\Omega\times(T,\infty))
    \quad\perogni T>0
 \end{equation}
 and the\/ {\rm Liapounov condition}
 \beeq{intLiap}
  \calG(u(t))
   \le \calG(u(0))
   \quext{for all~}\/t\ge 0.  
\end{equation}
\edde
\noindent%
Then, we have the following result, which will be proved in
the next Section~\ref{sec:utlimitata}:
\begin{teor}[Regularizing solutions]\label{teo:utlimitata}{\sl
 Let \eqref{alfadalbasso}, \eqref{W}--\eqref{coercW} and
 \eqref{regodata}--\eqref{regof} with \eqref{inizk} hold.
 Then, Problem~{\rm (P)} admits at least
 one regularizing solution. Moreover, there exist
 constants $\cc_1,\cc_2>0$ and a continuous and monotone
 function $\phi:[0,\infty)\to[0,\infty)$, all independent both 
 of the initial data and of time, and explicitly computable in terms of
 $\Omega$, $\alpha$, $W$, $f$, such that, for\/ {\rm every}
 regularizing solution and all $T>0$, it is
 \beal{regoui2}
   & \|u_t(t)\|^2_{L^\infty(\Omega)}\le 
    \cc_1\frac{1+\calG(u_0)}{T^{\cc_2}} \quad \perogni t\ge T,\\
  \label {regoui2.2}
   & d_\infty^2(u(t),0)\le
    \phi\Big(\cc_1\frac{1+\calG(u_0)}{T^{\cc_2}}\Big) \quad \perogni t\ge T.
 \end{align}
}
\end{teor}
\noindent%
In particular, thanks to the second inequality in \eqref{Liapd2}
and to \eqref{inizk}, the bounds \eqref{regoui2}, \eqref{regoui2.2}
depend only on the ``radius'' $\kappa$ of the initial datum with
respect to $d_2$.

Theorem~\ref{teo:utlimitata} is the starting point for
all the subsequent investigations. As a first
consequence, using the last of \eqref{regoui2} and \eqref{coercW},
from straightforward arguments there follows the
\begin{coro}[Separation]\label{prop:separation}{\sl
 Let \eqref{alfadalbasso}, \eqref{W}--\eqref{coercW} and
 \eqref{regodata}--\eqref{regof} hold, and let $u$ 
 be a regularizing solution. Then, for any $T>0$ there exist
 $\subr<0,\barr>0$,
 with $\inf I<\subr<0<\barr<\sup I$, such that 
 \beeq{sepsu}
    \subr\le u(x,t)\le \barr\quad\perogni x\in \Omega,~t\ge T.
 \end{equation}
}
\end{coro}
\beos
 The {\sl separation property}\/ \eqref{sepsu} stated in the Corollary
 improves the analogous inequality shown in \cite[Prop.~2.10]{SSS1} 
 and holding for less regular solutions (i.e., $\calX_2$-solutions in our
 notation) under additional assumptions on $W$.
\eddos
\noindent%
The local Lipschitz continuity of $W'$ (following from
\eqref{W}) and the simple argument used in 
\cite[Proof of Thm.~2.11]{SSS1} permit then to 
obtain immediately the 
\begin{coro}[Uniqueness]\label{unicita}{\sl
 Assume \eqref{alfadalbasso},
 \eqref{W}--\eqref{coercW} and
 \eqref{regodata}--\eqref{regof}. Let also
 $u,v$ be a pair of $\calX_2$-solutions
 corresponding 
 respectively to the initial data $u_0,v_0\in\calX_2$.
 Assume that, for some $s,c\ge 0$, 
 \beeq{prouniq}
   u(s)=v(s), \qquad
   d_\infty(u(t),0)+d_\infty(v(t),0)\le c
    \quad\perogni t\ge s,
 \end{equation}
 with $c$ independent of $t$. Then, $u\equiv v$
 on $[s,\infty)$.
}
\end{coro}
\noindent%
The proof of the next result will be detailed in 
Section~\ref{sec:utlimitata}.
\begin{coro}\label{semiflow}{\sl
 Under assumptions \eqref{alfadalbasso},
 \eqref{W}--\eqref{coercW} and
 \eqref{regodata}--\eqref{regof}, the set $\calS$ 
 of regularizing solutions to
 Problem~{\rm (P)} is a semiflow,
 whose space of regularized values is contained into
 $\calXi$.
}
\end{coro}
\beos
 Comparing our assumptions on $\alpha$, $W$ with those
 taken in \cite{SSS1}, we point out that here 
 (cf.~\eqref{coercW}), if $I\not=\RR$,  
 we are not able to consider
 potentials {\sl bounded}\/ in $\overline I$
 (like, e.g., the ``double obstacle''
 $W(r)\sim I_{[-1,1]}(r)-\lambda r^2/2$, $I_{[-1,1]}$ being
 the {\sl indicator function}\/ of $[-1,1]$). More precisely,
 this restriction is not required in the proof of 
 Theorem~\ref{teo:utlimitata}, where only 
 \eqref{W} is used, but in the subsequent 
 Corollaries~\ref{prop:separation} and~\ref{unicita}.
 Concerning $\alpha$, differently from 
 \cite{SSS1}, we cannot consider here the case in which 
 $\alpha$ is a maximal monotone function with
 some multivalued branch, and in particular we are not 
 able to deal with the situation where
 the domain of $\alpha$ is {\sl strictly}\/ included 
 in $\RR$ (as it happens, e.g., in the application to 
 {\sl irreversible}\/ phase transitions 
 considered in \cite{FV,LSS,LSS2}). Indeed, 
 in case $\dom \alpha\not=\RR$, 
 one can still deduce \eqref{regoui2}, but not 
 \eqref{regoui2.2} which is crucial for the long
 time analysis.
\eddos
\noindent%
\beos\label{ancorarego}
 The non-uniqueness of solutions to (P) 
 can be precised as follows.
 Given an initial datum $u_0\in\calX_2$, from it more than
 one solution can emanate. In particular, there are
 one, or more, regularizing solutions starting from $u_0$, 
 at least one of which is limiting, and all these 
 regularizing solutions are taken as elements of the semiflow $\calS$. 
 Other solutions can also exist which are not elements of $\calS$. 
 In particular, (nonlimiting) smooth solutions 
 enjoying \eqref{regoui2}
 but not \eqref{intLiap} are excluded from~$\calS$.
\eddos
\noindent%
Let us now come to the long time behavior.
\begin{teor}[Global attractor]\label{teoattrattore}{\sl
 Assume \eqref{alfadalbasso}, \eqref{W}--\eqref{coercW} and
 \eqref{regodata}--\eqref{regof}. Then,
 the semiflow $\calS$ associated 
 with Problem~{\rm (P)} admits the
 global attractor $\mathcal{A}$, which is compact 
 in $\calX_2$ and ``sequentially weakly compact'' in 
 $\calXi$ (i.e., sequences in $\mathcal{A}$ admit
 subsequences ``weakly'' converging in $\calXi$).
 }
\end{teor}
\begin{teor}[Exponential attractors]\label{teoattrattoreesp}{\sl
 Suppose that \eqref{alfadalbasso},
 \eqref{W}--\eqref{coercW} and
 \eqref{regodata}--\eqref{regof} hold.
 Then, the semiflow $\calS$ associated 
 with Problem~{\rm (P)} admits an exponential
 attractor $\mathcal{M}$. More precisely, $\calM$ is 
 a compact subset of $V$ which 
 attract exponentially fast with respect to the
 $V$-norm any $d_2$-bounded subsets of $\calX_2$.
 }
\end{teor}
\noindent%
\beos\label{VnonX}
 We decided to show existence of $\calM$ by working 
 in $V$ rather than in $\calX_2$ just for the 
 sake of simplicity. Indeed, reinforcing the 
 differentiability assumptions on $\alpha$ and $W$
 and refining the estimates in Section~\ref{secexpo} 
 (roughly speaking, we could put one more Laplacian 
 in the test functions used there), 
 it should be possible to obtain estimates analogous
 to \eqref{contod19}--\eqref{contod20} below, but with one 
 more order of space derivatives inside. We omitted
 to perform such a procedure since it would
 be rather lenghty and technical. 
 Its outcome would be the existence of
 an exponential attractor in $H^2(\Omega)$
 and, in fact, in $\calX_2$ (indeed, the
 contribution of $W'$ in $d_2$ is automatically
 controlled by uniform eventual $\calXi$-boundedness
 of solutions and local Lipschitz continuity of $W'$).
\eddos
\noindent%
As recalled in Section~\ref{preli},
the existence of $\calM$ entails that the global attractor
$\calA$ is contained in $\calM$ and has finite fractal
dimension in $V$ (actually in $H^2(\Omega)$ on account
of the Remark above).

As a final issue, by virtue of the $L^\infty$-bound 
on $u_t$, we are able to sharpen
the results in \cite{SSS1} concerning $\omega$-limits of the elements
of $\calS$. Actually, since $\alpha(0)=0$, it is clear
(cf.~\cite[Thm.~2.13]{SSS1}) that the stationary 
states $u_\infty$ of \eqref{eqn} are solutions of 
\begin{equation}\label{eq-stazionaria}
  B u_\infty+W'(u_\infty)=f \quext{ in }\,H.
\end{equation}
It is well known that, since $W$ needs not be
convex, \eqref{eq-stazionaria} may well admit infinitely
many solutions \cite{haraux}, all of which, due
to \eqref{W}, \eqref{coercW} and standard elliptic regularity
results, belong to $\calXi$. Thus, given $u\in \calS$, the question of the 
convergence of all the trajectory $u(t)$ to one of these solutions may be non
trivial. As in \cite{SSS1}, we are able to show this property by making use 
of the so-called {\L}ojasiewicz-Simon inequality \cite{Lo1,Lo2,Si2}, 
at least provided that
\beeq{analit}
   W|_{I_0} \quext{is\/ {\rm real analytic}},
\end{equation}
where $I_0\subset I$ is an open interval containing $0$ and
such that $W'(r)r>0$ for all $r\in I\setminus I_0$. Clearly,
$I_0$ exists thanks to \eqref{coercW}; moreover, by maximum principle 
arguments, any solution to \eqref{eq-stazionaria} takes values
in a compact subset of $I_0$. Then, we have the following
\begin{teor}[Convergence to the stationary states]\label{Loj}{\sl
 Let us assume hypotheses\/ \eqref{alfadalbasso},
 \eqref{W}--\eqref{coercW}, \eqref{regodata}--\eqref{regof}
 and\/ \eqref{analit}.
 Then, letting $u$ be a regularizing solution, the $\omega$-limit 
 of\/ $u$ consists of a\/ {\rm unique}
 function $\ui$ solving \eqref{eq-stazionaria}.
 Furthermore, as $t\nearrow+\infty$,
 \beeq{coinfty2}
    u(t) \to \ui \qquext{strongly in }\,V\cap C(\barO),
 \end{equation}
 i.e., we have convergence 
 for the\/ {\rm whole trajectory $u(t)$}.
 }
\end{teor}
\noindent%
The difference between this result and 
\cite[Thm.~2.18]{SSS1} lies
in the fact that, thanks to \eqref{regoui}, we need not assume any 
growth condition on $\alpha$. Roughly speaking, the $L^\infty$-bound on $u_t$ combined
with the regularity and the coercivity of $\alpha$ (see \eqref{alfadalbasso}) reduces the 
nonlinearity $\alpha$ to an almost ``linear'' contribution
and makes the analysis of the convergence
of the trajectory simpler. In fact, Theorem \ref{Loj} can be 
proved by simply adapting the proof given in~\cite{CJ}. We leave the details
to the reader.
%
\begin{osse}[The asymptotically autonomous case]\label{nonauto}\rm
 For the sake of studying $\omega$-limits, we could also consider  
 time dependent sources, by assuming, instead of \eqref{regof},
 \begin{equation}\label{nonaut1}
   f\in L^2(0,+\infty;L^{\infty}(\Omega)), \quad
    f_t\in L^1(0,+\infty;L^\infty(\Omega)).
 \end{equation}
 Indeed, it could be shown that Theorem~\ref{teo:utlimitata} 
 and Corollaries~\ref{prop:separation},~\ref{unicita},
 and~\ref{semiflow} still hold in this setting.
 Moreover, assuming also that there exist $c,\xi>0$ such that 
 \beeq{decadef}
   t^{1+\xi}\int_t^\infty \|f(s)\|^2\,\dis \le c 
    \quext{for all }\,t\ge0,
 \end{equation}
 Theorem~\ref{Loj} could be extended as well (see also
 \cite{CJ,GPS2} for this kind of assumptions).
\end{osse}


\section{Regularization for strictly positive times}
 \label{sec:utlimitata}

{\bf Proof of Theorem~\ref{teo:utlimitata}}.~~%
We shall use an Alikakos-Moser \cite{Al} iteration argument
 for which some
a priori estimates are needed. In particular, we shall work on the 
(formal) time derivative of \eqref{eqn}, given by
\beeq{eqnprimo}
  \alpha'(u_t)u_{tt}+B u_t+W''(u)u_t=0.
\end{equation}
Of course, \eqref{eqnprimo} needs not make sense
if $u$ is just an $\calX_2$-solution. However, 
we can write it for Problem~(P$_n$), derive
the estimates at the level $n$, and then let them
pass to the limit $n\nearrow\infty$ 
using the semicontinuity properties of norms
w.r.t.~weak convergences.  This approach
has the drawback that, at a first stage,
the estimates will hold only for the ``limiting
solutions''. They will be properly extended 
to all regularizing solutions in the second
part of the proof.

Before proceeding, we introduce some further
notation. For simplicity, we 
shall omit the index $n$ of the approximation
in all what follows. The symbol $c$ will stand for a positive
constant, possibly varying even inside one single 
line, which is allowed to depend on the data
$\Omega$, $\alpha$, $W$, $f$, but neither on the initial
values, nor on time. The constant(s) $c$ will be also 
independent of the exponents $p_j$
of the iteration process (see below) and, of 
course, of $n$. Some $c$'s whose precise value is needed 
will be distinguished by 
noting them as $c_i$, $i\ge 0$. Let us now set,
for $p\in[2,\infty)$,
\beeq{defiap}
  a_p(s):=\int_0^s\alpha'(r)|r|^{p-2}r\,\dir
\end{equation}
and notice that (recall that  $\alpha(0)=0$)
\beeq{apsu}
  \frac\sigma{p}|s|^p\le a_p(s)\le\alpha(s)|s|^{p-2}s \quad\perogni s\in\RR.
\end{equation}
Moreover, it is clear that (at least formally, as noted above)
\beeq{ap2}
  \ddt a_p(u_t)=\alpha'(u_t)|u_t|^{p-2}u_t\+u_{tt}.
\end{equation}
Then, testing \eqref{eqnprimo} by $u_{t}$, recalling the second of
\eqref{W} and adding $\lambda\|u_t\|^2$ on both hands sides, 
and integrating over $(0,t)$, we get
\beeq{regoV-0}
  2\|a_2(u_t(t))\|_{L^1(\Omega)}
   +2\|u_t\|_{L^2(0,t;V)}^2
  \le 2\|a_2(u_t(0))\|_{L^1(\Omega)}
   +c\|u_t\|_{L^2(0,t;H)}^2.
\end{equation}
To control the latter term in the \rhs\ above, we can use 
\eqref{stimaL2H}. The other one,
by \eqref{apsu} with $p=2$ and Young's inequality,
becomes
\beeq{ut0}
  2\|a_2(u_t(0))\|_{L^1(\Omega))}
   \le \|\alpha(u_t(0))\|^2+\|u_t(0)\|^2
   \le c(1+\kappa)^2,
\end{equation}
where the latter inequality is a consequence of a comparison
in \eqref{eqn} (written for \Pn)
and of assumption \eqref{regodata} ($\kappa$ is as in 
\eqref{inizk}). Actually, $\alpha^{-1}$
is Lipschitz continuous due to \eqref{alfadalbasso}.
In conclusion, from \eqref{regoV-0} we obtain
\beeq{regoV}
  2\|a_2(u_t)\|_{L^\infty(0,\infty;L^1(\Omega))}
   +2\|u_t\|_{L^2(0,\infty;V)}^2\le c_0(1+\kappa)^2.
\end{equation}
We can now describe the two estimates which are at the core of 
the iteration process. 

\vspace{2mm}

\noindent%
{\bf First estimate.}~~%
Let $j\ge 1$, $p_j>1$, and let us test \eqref{eqnprimo}
by $|u_t|^{p_j-2}u_t$, so that 
\beeq{conto11}
  \ddt \io a_{p_j}(u_t)
   +\big(Bu_t,|u_t|^{p_j-2}u_t\big)
   \le \lambda\|u_t\|_{p_j}^{p_j}
\end{equation}
(we agree, here and in the sequel, to note by
$\|\cdot\|_p$ the norm in $L^p(\Omega)$
for $p\in[1,\infty]$). By definition of $B$
and Poincar\'e's inequality (everything works with minor
changes also in the Neumann case),
\beeq{conto12}
  \big(Bu_t,|u_t|^{p_j-2}u_t\big)
   \ge 
   \frac{4(p_j-1)}{p_j^2}\io \Big|\nabla\big( |u_t|^{\frac{p_j-2}2}u_t\big)\Big|^2
   \ge \frac{c_1}{p_j}\|u_t\|_{3p_j}^{p_j},
\end{equation}
for some $c_1>0$.
Assuming then that there exist $T_j,\ell_j>0$ such that 
\beeq{conto13} 
  p_j\|a_{p_j}(u_t)\|_{\Tuno(T_j,\infty;L^1(\Omega))}\le \ell_j, \qquad
  p_j\|u_t\|_{\Tpj(T_j,\infty;L^{p_j}(\Omega))}^{p_j}\le \ell_j
\end{equation}
%
and multiplying \eqref{conto11} by $p_j$,
from Lemma~\ref{unigronw} we get, for $\tau_j\in(0,1]$,
\beeq{st11}
  p_j\|a_{p_j}(u_t(t+\tau_j))\|_{L^1(\Omega)}
   \le \ell_j\Big(\lambda+\frac1{\tau_j}\Big)
   \qquad\perogni t\ge T_j,
\end{equation}
whence, recalling \eqref{apsu}, we also have 
\beeq{st11bis}
  \|u_t(t+\tau_j)\|_{p_j}^{p_j}
   \le \frac{\ell_j}\sigma\Big(\lambda+\frac1{\tau_j}\Big)
   \qquad\perogni t\ge T_j.
\end{equation}
Moreover, integrating $p_j$ times \eqref{conto11}
over $(t,t+1)$ for $t\ge T_j+\tau_j$, and taking 
\eqref{conto12}, \eqref{st11} 
into account, it is not difficult to infer
\beeq{st12}
  \int_t^{t+1} \|u_t(s)\|_{3p_j}^{p_j}\,\dis
   \le \frac{\ell_j}{c_1}\Big(2\lambda+\frac1{\tau_j}\Big)
   \qquad\perogni t\ge T_j+\tau_j.
\end{equation}

\vspace{2mm}

\noindent%
{\bf Interpolation argument.}~~%
By elementary interpolation of $L^p$ spaces, we have
\beeq{int1}
  \|u_t(t)\|_{7p_j/3}
   \le \|u_t(t)\|_{p_j}^{1/7}\|u_t(t)\|_{3p_j}^{6/7}
   \qquad\perogni t\ge T_j+\tau_j.
\end{equation}
Hence, still for $t\ge T_j+\tau_j$,
\beeq{int2}
  \int_t^{t+1}\|u_t(s)\|_{7p_j/3}^{7p_j/6}\,\dis
   \le \|u_t\|_{L^\infty(t,t+1,L^{p_j}(\Omega))}^{p_j/6}
   \int_t^{t+1}\|u_t(s)\|_{3p_j}^{p_j}\,\dis.
\end{equation}
Thus, from \eqref{st11bis} and \eqref{st12},
\beeq{int3}
  \|u_t\|_{\calT^{7p_j/6}(T_j+\tau_j,\infty;L^{7p_j/3}(\Omega))}^{7p_j/6}
   \le \Big(\frac{\ell_j}{\sigma}\Big)^{1/6}
       \Big(\lambda+\frac1{\tau_j}\Big)^{1/6}
     \frac{\ell_j}{c_1}\Big(2\lambda+\frac1{\tau_j}\Big).
\end{equation}
In conclusion, there exists $c_2$ depending only on $c_1,\sigma,\lambda$ 
and such that
\beeq{int4}
  \|u_t\|_{\calT^{7p_j/6}(T_j+\tau_j,\infty;L^{7p_j/3}(\Omega))}^{p_j}
   \le c_2\ell_j\Big(1+\frac1{\tau_j}\Big).
\end{equation}

\vspace{2mm}

\noindent%
{\bf Second estimate.}~~%
We now test \eqref{eqn} by $|u_t|^{q-2}u_t$, with $q>1$ to be chosen
later. Owing to the bound \eqref{regou2} and using \eqref{alfadalbasso},
it is clear that
\beeq{conto21}
  \io\alpha(u_t)|u_t|^{q-2}u_t
   \le \|-Bu-W'(u)+f\|_2\+\|u_t\|_{2q-2}^{q-1}
   \le c(1+\kappa)\|u_t\|_{2q-2}^{q-1}.
\end{equation}
Consequently,
\beeq{conto21.2}
 \sigma\|u_t\|_q^q
   \le c(1+\kappa)\|u_t\|_{2q-2}^{q-1}.
\end{equation}
The above relations \eqref{conto21}--\eqref{conto21.2}
hold pointwise in $t$.
Then, integrating \eqref{conto21} over $(t,t+1)$ for 
$t$ greater than a suitable $S$ and using
the latter inequality in \eqref{apsu}, 
we get, for some $c_3$ depending only on $C,\sigma$,
\beeq{conto22} 
  q\|a_q(u_t)\|_{\Tuno(S,\infty;L^1(\Omega))}
   +q\|u_t\|_{\Tq(S,\infty;L^q(\Omega))}^q
  \le c_3 q(1+\kappa)\int_t^{t+1}\|u_t(s)\|_{2q-2}^{q-1}\,\dis.
\end{equation}

\vspace{2mm}

\noindent%
{\bf Bootstrap.}~~%
At this point, if we take in the previous argument
\beeq{sceltaq}
  S=T_{j+1}:=T_j+\tau_j, \qquad
   q=p_{j+1}:=\frac{7p_j}6+1,
\end{equation}
relation \eqref{conto22} is readily rewritten as 
\bealo
  & p_{j+1}\|a_{p_{j+1}}(u_t)\|_{\Tuno(T_{j+1},\infty;L^1(\Omega))}
   +p_{j+1}\|u_t\|_{\calT^{p_{j+1}}(T_{j+1},\infty;L^{p_{j+1}}(\Omega))}^{p_{j+1}}\\
 \label{conto22.2} 
  & \mbox{}~~~~~~~~~~
    \le c_3 p_{j+1}(1+\kappa)\int_t^{t+1}\|u_t(s)\|_{2p_{j+1}-2}^{p_{j+1}-1}\,\dis.
\end{align}
Hence, recalling \eqref{int4}, the \lhs\ above is majorized by
\beeq{conto22.3} 
   c_3 p_{j+1} (1+\kappa)c_2^{7/6}\ell_j^{7/6}\Big(1+\frac1{\tau_j}\Big)^{7/6}
     \le c_4\ell_j^{7/6}p_j\Big(1+\frac1{\tau_j}\Big)^{7/6}(1+\kappa).
\end{equation}
Thus, we can define
\beeq{ellpiu1}
  \ell_{j+1}:=c_4\ell_j^{7/6}p_j\Big(1+\frac1{\tau_j}\Big)^{7/6}(1+\kappa),
\end{equation}
so that \eqref{conto22.3} implies \eqref{conto13} at the step $j+1$.
More precisely, since by \eqref{regoV} we can take 
\beeq{primivalori}
   T_1:=0, \qquad p_1:=2, \qquad \ell_1:=c_0(1+\kappa)^2,
\end{equation}
assuming that $\epsilon\in(0,1)$ is given, we also choose
\beeq{defitauj}
   \tau_j:=\frac\epsilon{j^2}, \quext{so that }\,
   T_{j+1}=T_{j}+\tau_j \le c\epsilon
   \quad\perogni j\ge1
\end{equation} 
and for $c>0$ independent of~$j$.
At this point, let us set, for notational simplicity,
\beeq{newnot}
  b:=7/6, \qquad
   B_j:=\sum_{i=0}^jb^i\le 6b^{j+1}.
\end{equation} 
Then, it is not difficult to get from \eqref{ellpiu1} 
(cf.~also \eqref{primivalori})
\beeq{ellpiu1.2}
  \ell_{j+1}\le c_4^{B_{j-1}}c_0^{b^j}(1+\kappa)^{p_{j+1}}
   \prod_{i=1}^j p_i^{b^{j-i}}
   \prod_{i=1}^j \Big(1+\frac{i^2}{\epsilon}\Big)^{b^{j-i+1}},
\end{equation}
whence, noting that 
\beeq{maggiorp}
   c_5 b^j\le p_j\le c_6 b^{2j} \quad\perogni j\ge 1
\end{equation}
and for some $c_5,c_6>0$ independent of $j$,
and passing to the 
logarithm, it is not difficult to show that
\bega{stserie1}
  \Big(\prod_{i=1}^j p_i^{b^{j-i}}\Big)^{1/p_j}
   \le c,\\
\label{stserie2}
  \bigg(\prod_{i=1}^j
  \Big(1+\frac{i^2}{\epsilon}\Big)^{b^{j-i+1}}\bigg)^{1/p_j}
   \le \frac{c}{\epsilon^{c_7}}.
\end{gather}
%
%
Collecting the above estimates, we infer
\beeq{ellpiu1.3}
  \ell_{j+1}^{1/p_{j+1}}\le \frac{c(1+\kappa)}{\epsilon^{c_7}}.
\end{equation}
Thus, 
\eqref{st11bis} (written at the step $j+1$)
gives, for all $j\in \NN$,
\beeq{key1}
  \|u_t(t)\|_{p_j} \le \frac{c(1+\kappa)}{\epsilon^{c_8}}
   \qquad\perogni t\ge T_{j+1}.
\end{equation}
From \eqref{int4} we also have
\beeq{key2}
  \|u_t\|_{\calT^{p_{j+1}-1}(T_{j+1},\infty;L^{2(p_{j+1}-1)}(\Omega))}
   \le \frac{c(1+\kappa)}{\epsilon^{c_8}}.
\end{equation}
Finally, taking the limit of \eqref{key1} as $j\nearrow\infty$
we obtain
\beeq{key3}
  \|u_t(t)\|_\infty
   \le  \frac{c_9(1+\kappa)}{\epsilon^{c_8}}
    \quad\perogni t\ge c\epsilon,
\end{equation}
where the last $c$ is the same as in \eqref{defitauj}.
Hence, by arbitrariness of $\epsilon$,
$u_t(t)$ is essentially bounded for a.e.~$t>0$. More precisely,
squaring \eqref{key3}, recalling \eqref{inizk},
and owing also to the first inequality
in \eqref{Liapd2}, \eqref{regoui2} follows at once.
Recalling \eqref{alfadalbasso}, and using in particular
that $\alpha$ is defined on the whole real line, 
we also obtain
\beeq{regoprepa}
  \|\alpha(u_t)\|_\infty
   \le\phi\Big(\cc_1\frac{1+\calG(u_0)}{T^{\cc_2}}\Big) 
   \quad \perogni t\ge T,
\end{equation}
where $\phi$ depends only on $\alpha$.
Then, rewriting \eqref{eqn} as 
\beeq{elli}
   Bu+W'(u)+\lambda u=f+\lambda u-\alpha(u_t),
\end{equation}
and viewing it as a time dependent family of 
elliptic problems with monotone nonlinearity and
uniformly bounded forcing term, it is not difficult
to obtain also \eqref{regoui2.2} as a consequence 
of standard maximum principle arguments.
More precisely, one can test \eqref{elli} by
$|W'(u)+\lambda u|^{p-2}(W'(u)+\lambda u)$ for
$p\in[2,\infty)$ and then let $p\nearrow\infty$.

To conclude the proof of Theorem~\ref{teo:utlimitata}, 
we recall that the procedure above has to be intended in 
the framework of Problem~(P$_n$). Then, the bounds
\eqref{regoui2}, \eqref{regoui2.2}, as well as the 
Liapounov condition \eqref{intLiap}, pass easily
to the limit $n\nearrow\infty$ thanks to 
lower semicontinuity of norms with respect to 
weak and weak star convergences.
More precisely, to obtain \eqref{intLiap} the 
following property (of straightforward proof) 
is used:
\bele\label{lemmalsc}
 The functional $\calG$ is sequentially 
 weakly lower semicontinuous in $\calX_2$,
 namely we have
 \beeq{lsc1}
   \calG(u) \le \liminf_{n\nearrow\infty}\calG(u_n)
 \end{equation}
 if $\{u_n\}\subset\calX_2$ tends to some limit $u$
 {\rm weakly} in $\calX_2$. The same property
 holds also for $\calF$.
\enle
\noindent%
The proof of Theorem~\ref{teo:utlimitata} is however not
yet complete since, up to now, we have just showed that 
any {\sl limiting solution}\/ is a 
{\sl regularizing solution} and fulfills \eqref{regoui2},
\eqref{regoui2.2} and \eqref{intLiap}. To conclude,
we have to prove that {\sl any}\/ regularizing
solution $u$ (i.e.~also a {\sl nonlimiting} one) satisfies
\eqref{regoui2} and \eqref{regoui2.2}
(while \eqref{intLiap} is now postulated
in Definition~\ref{defisolXi}).
Here, the key point is to notice that, by \eqref{regoui}
and Cor.~\ref{unicita}, taken any $s>0$, from the ``datum''
$u(s)$ at most one solution emanates. Thus, 
any regularizing $u$ is also
``limiting'' as it is restricted to $[s,\infty)$.
This means that, referring for instance to \eqref{regoui2},  
we have at least
\beeq{regoui2st}
  \|u_t(t)\|^2_{L^\infty(\Omega)}\le 
    \cc_1\frac{1+\calG(u(s))}{(T-s)^{\cc_2}}
    \quad\perogni t\ge T>s>0.
\end{equation}
Then, \eqref{regoui2} follows easily by first using
\eqref{intLiap} (with $s$ in place of $t$) and then 
taking the limit for $s\searrow 0$. The bound
\eqref{regoui2.2} is proved exactly in the same way 
and concludes the proof of 
Theorem~\ref{teo:utlimitata}.\dimbox
\beos\label{Liapaltritempi}
 Notice that, for {\sl any}\/ regularizing solution, 
 there holds the property (slightly stronger
 than \eqref{intLiap})
 \beeq{intLiapst}
  \calG(u(t))
   \le \calG(u(s))
   \quext{for all~}\/t\ge s\ge 0.
 \end{equation}
 Indeed, if $s=0$, then \eqref{intLiapst}
 reduces to \eqref{intLiap}. Otherwise,
 $u$ coincides on $[s,\infty)$ with a
 limiting solution. Thus, \eqref{intLiapst} can 
 be shown by noting as before that  
 $u$ is limiting on $[s,\infty)$, 
 considering (P$_n$) w.r.t.~the ``initial''
 datum $u(s)$, and finally letting 
 $n\nearrow\infty$.
\eddos
\noindent%
{\bf Proof of Corollary~\ref{semiflow}.}~~%
Property {\sl (S1)}\/ is evident and {\sl (S4)}\/ follows
from Cor.~\ref{unicita}. Next, {\sl (S2)}\/ and
{\sl (S3)}\/ are immediate once one notes that 
$v$ (in {\sl (S2)}\/) and $z$ (in {\sl (S3)}\/)
fulfill \eqref{intLiap} thanks to 
Remark~\ref{Liapaltritempi}. Finally, let us prove 
{\sl (S5)}. Although we could use here the regularization
properties \eqref{regoui2}, \eqref{regoui2.2}, 
we rather give a proof which essentially
relies only on \eqref{regou2},
since we think it is interesting to notice that
the strong-weak semicontinuity properties require
no smoothing effect.

Thus, to show the first of {\sl (S5)}, 
we start by observing that, due to
\eqref{regou}, any $u\in\calS$ stays in
$C_w([0,\infty);H^2(\Omega))$, so that
we just have to prove that, as $s,t\in [0,\infty)$
and $s$ tends to $t$, $(W'+\lambda)(u(s))$
goes to $(W'+\lambda)(u(t))$ weakly in $H$.
To see this, we first notice (cf.~also~\cite[Sec.~6]{RS})
that there exists $c\ge 0$ such that 
$\|(W'+\lambda)(u(s))\|\le c$ for {\sl all}\/
(not just a.e.) $s\in[0,\infty)$.
Then, it is clear that, as $s\to t$, any subsequence 
of $(W'+\lambda)(u(s))$
admits a subsequence weakly convergent in $H$,
whose limit is identified as $(W'+\lambda)(u(t))$
thanks to the convergence $u(s)\to u(t)$,
which is strong in $H$, the monotonicity
of $W'+\lambda\Id$, and \cite[Lemma~1.3, p.~42]{Ba}.
This proves weak continuity of single trajectories.
Note that if we admit use of \eqref{regoui2}, \eqref{regoui2.2}, 
we actually get more, namely $W'(u(\cdot))$ 
is strongly continuous with values
in $C(\barO)$ at least for strictly positive times.

To conclude, let us show the second property in {\sl (S5)}.
Letting then $u_n,u\zzn$ as in {\sl (S5)}, as $u\zzn$
tends to $u_0$ in $\calX_2$, it is in particular 
bounded in $\calX_2$. This entails that 
\eqref{regou2}, \eqref{regoui2}, \eqref{regoui2.2}
hold uniformly in $n$. By compactness arguments
(similar to those in \cite[Subsec.~3.3]{SSS1})
and using \cite[Cor.~4]{Si},
we then obtain that (a not relabelled subsequence of) 
$u_n$ satisfies, for all $T>0$,
\beal{cou1}
  & u_n\to u \quext{strongly in }\,\CZV,\\
 \label{cou2}
  & (W'+\lambda)(u_n)\to (W'+\lambda)(u)
   \quext{weakly in }\,\LDH,
\end{align}
where $u$ is an $\calX_2$-solution to
Problem~(P) with initial datum $u_0$, and it
satisfies \eqref{regou2}, \eqref{regoui2}
and \eqref{regoui2.2}. In particular, given
any $t>0$, by \eqref{cou1} $u_n(t)$ goes to 
$u(t)$ strongly in $V$. Then, by uniform boundedness, 
this convergence is also weak in $H^2(\Omega)$.
As before, the monotonicity of $W'+\lambda\Id$
and the bound $\|(W'+\lambda)(u_n(t))\|\le c$, which
is uniform both in $n$ and in $t$, permit to show that 
$(W'+\lambda)(u_n(t))\to (W'+\lambda)(u(t))$ weakly in
$H$ (no further extraction of subsequence is required
here, since the limit is already identified).
To conclude, we have to see that $u$ is a regularizing
solution (i.e.~it also fulfills \eqref{intLiap}). To
prove this, it suffices to write \eqref{intLiap}
for $u_n$ and take the liminf as $n\nearrow\infty$. 
Indeed, the left hand 
side can be treated by Lemma~\ref{lemmalsc},
while the \rhs\ passes directly to the limit since 
$u\zzn\to u_0$ {\sl strongly}\/ in $\calX_2$ and
it is easy to check that $\calG$ is continuous
with respect to $d_2$.\dimbox


\section{Long time behavior}
\label{esiattra}

{\bf Proof of Theorem~\ref{teoattrattore}.}~~%
We shall show the following facts:\\[2mm]
{\sl (L1)}~~The semiflow $\calS$ possesses a 
   Liapounov function;\\[1mm]
{\sl (L2)}~~The set of stationary points of $\calS$ is bounded
  in $\calX_2$;\\[1mm]
{\sl (L3)}~~The semiflow $\calS$ is {\sl asymptotically compact}, namely
 for any sequence $\left\{u_n^{0}\right\}_{n\in\mathbb{N}}$ bounded
 in $\mathcal{X}_2$ and any positive sequence $\left\{t_n\right\}_{n\in\mathbb{N}}$,
 $t_n\nearrow \infty$, any sequence of the form $\{u_n(t_n)\}$, 
 where $u_n\in \calS$ and $u_n(0)=u_n^0$,
 is precompact in $\mathcal{X}_2$.\\[2mm]
By the theory of global attractors
(see, e.g., \cite[Theorem 3.2]{lady} or \cite[Thm.~5.1]{Ba1}),
{\sl (L1)--(L3)}\/ would imply the existence of a global
attractor compact in $\calX_2$. However, here neither
the ``standard'' theory in \cite{lady}, nor the 
``generalized''  theory in \cite{Ba1}, can be directly
applied since we have no uniqueness and
just strong-weak semicontinuity. Nevertheless, we shall
show in the Appendix that the validity of~\cite[Thm.~5.1]{Ba1}
can be extended also to this case.
\beos\label{ossedissi}
 The use of this method permits to bypass the direct proof
 of existence of an $\calX_2$-bounded absorbing set,
 which seems difficult to get here due to the possibly fast
 growth of $\alpha$ at $\infty$. Of course, 
 {\sl a posteriori}\/ the dissipativity
 property will be satisfied just as a consequence of
 the existence of the global attractor.
\eddos
\noindent%
To proceed, we first notice that, by the {\sl energy estimate}\/
(obtained testing \eqref{eqn} by $u_t$), $\calE$ is a 
Liapounov functional. Note that the regularity of {\sl any}\/
$\calX_2$-solution is sufficient to justify this 
estimate (and this is the reason why we do not use 
here the functional $\calG$, which also enjoys
a Liapounov property, at least for {\sl regularizing}\/
solutions, by Remark~\ref{Liapaltritempi}).
Thus, {\sl (L1)}\/ holds.
Second, {\sl (L2)}\/ is an easy consequence of
well-known elliptic regularity results, 
so that it just remains to show {\sl (L3)},
whose proof will be split in a number of steps.
\bele\label{step1}
 Given $0<\tau<T<\infty$, there exists $c$ depending
 on $\tau,T$ and on the initial datum
 such that any regularizing solution $u$
 satisfies the further bounds
 \beal{stnew11}
   & \|u_{tt}\|_{L^{2}(\tau,T;H)}+\|u_{t}\|_{L^{\infty}(\tau,T;V)}\le c,\\
  \label{stnew12}
   & \|Bu_{t}\|_{L^2(\tau,T;H)}\le c.
\end{align}
\enle
\begin{proof}
 We can prove \eqref{stnew11}--\eqref{stnew12} 
 by working on (P$_n$) and then
 letting $n\nearrow\infty$. As before, we omit the subscript
 $n$, for simplicity. Indeed, since we just 
 consider strictly positive times, $u$ can be thought
 as a limiting solution. In this regard,
 \eqref{stnew11} is obtained by testing \eqref{eqnprimo} 
 by $(t-\tau) u_{tt}$ 
 and using monotonicity of $\alpha$
 together with \eqref{stimaL2H} and \eqref{regoui}.
 Next, \eqref{stnew12} follows by making a comparison
 in \eqref{eqnprimo} and using \eqref{regoui} and 
 the uniform boundedness of $\alpha'(u_t)$ and
 $W''(u)$. The technical details 
 of the procedure, as well as the standard 
 argument for passing to the limit with $n$,
 are left to the reader.
\end{proof}
To proceed, we set, just to avoid some technicalities, 
$f\equiv 0$. We have the 
\bele\label{step2}
 Let $z\in\calS$. Setting, for $s>0$, 
 \beeq{newH}
   H(z(s)):=-\big(\alpha(z_t(s)),(Bz_t+W''(z)z_t)(s)\big)
    -\frac12\big(\alpha(z_t(s)),(Bz+W'(z))(s)\big),
 \end{equation}
 for any $\tau,M>0$ there holds
 \begin{equation}\label{as-comp2new}
   \mathcal{F}(z(\tau+M))=e^{-M}\mathcal{F}(z(\tau))
    +\int_{\tau}^{\tau+M}e^{s-\tau-M}H(z(s))\,\dis.
 \end{equation}
\enle
\begin{proof}
 Since we work on $[\tau,\infty)$, we can use the 
 further regularity properties 
 \eqref{stnew11}--\eqref{stnew12},
 which allow us to test \eqref{eqn} by 
 $(Bz_t+W''(z)z_t)+\frac{1}{2}(Bz+W'(z))$.
 Integrating over $(\tau,\tau+M)$, 
 we readily get~\eqref{as-comp2new}.
\end{proof}
\beos
 Let us note that, using, e.g., 
 \cite[Lemme~3.3, p.~73]{Br}, we 
 get more precisely that the function
 $t\mapsto\calF(z(t))$ is absolutely continuous
 on $[\tau,\infty)$ for all $\tau>0$. This permits,
 in particular, to improve (in our specific case)
 the first condition in {\sl (S5)}. Namely,
 the elements of our semiflow $\calS$ 
 belong to $C((0,\infty);\calX_2)$ 
 (compare this fact with condition \cite[(C1)]{Ba1}).
\eddos
\noindent%
Let us now complete the proof of {\sl (L3)}. We use here the
``energy method'' originally devised by Ball in \cite{Ba2}
(see also~\cite{moise-rosa-wang04} for an 
extension to nonautonomous systems). 
Take $\tau,M$ as before, and let
$v_n$ be the (unique) regularizing solution satisfying,
for $t\in[0,\infty)$,
$v_n(t)=u_n(t_n+t-M-\tau)$ (so that, in particular,
$v_n(0)=u_n(t_n-M-\tau)$, $v_n(\tau)=u_n(t_n-M)$ and
$v_n(\tau+M)=u_n(t_n)$). Since by \eqref{regoui2.2} there
exists $k>0$ such that $d_\infty(v_n(t),0)\le k$ for all
$n\in\NN$ and $t\in[0,\infty)$, by weak compactness
we have that there exist $\chi_{-M},\chi\in \calX_2$ such
that $v_n(\tau)\to \chi_{-M}$ and $v_n(\tau+M)\to \chi$
weakly in $\calXi$. Then, writing \eqref{as-comp2new} 
for $z=v_n$, we get
\begin{align}\nonumber
  \mathcal{F}(u_n(t_n))-e^{-M}\mathcal{F}(u_n(t_n-M))
   & = \mathcal{F}(v_n(\tau+M))-e^{-M}\mathcal{F}(v_n(\tau))\\
 \label{as-comp3new}
  & = \int_{\tau}^{\tau+M}e^{s-\tau-M}H(v_n(s))\,\dis
   =:\calH(v_n).
\end{align}
Next, let us note that, at least up to a not relabelled
subsequence, $v_n$ properly tends to an $\calX_2$-solution $v$.
Thus, in particular, we have that $v(\tau)=\chi_{-M}$ and 
$v(\tau+M)=\chi$.
Moreover, still by \eqref{regoui2.2}, 
$d_\infty(v(t),0)\le k$ for all $t\in[0,\infty)$.
Thus, setting $v_0:=\lim_{n\nearrow\infty}v_n(0)$, 
since by the existence property there exists at least
one $\zeta\in \calS$ such that $z(0)=v_0$, by 
Corollary~\ref{unicita} it must be $\zeta\equiv v$ on $[0,\infty)$,
which means that $v$ is itself an element of $\calS$ and,
consequently, satisfies \eqref{as-comp2new}.
Thus, noting that, by \eqref{stnew11}, \eqref{stnew12}
and weak compactness, $\calH(v_n)$ tends 
to $\calH(v)$, taking the $\limsup$ 
in~\eqref{as-comp3new} one gets
\begin{align}\nonumber
  \limsup_{n\nearrow\infty}\mathcal{F}(u_n(t_n))
   & \le c e^{-M} + \limsup_{n\nearrow\infty}\mathcal{H}(v_n)\\
 \nonumber
  & = c e^{-M} + \mathcal{H}(v)\\
 \nonumber
  & = c e^{-M} + \calF(v(\tau+M))-\calF(v(\tau))e^{-M}\\
 \label{as-comp4new}
  & \le c e^{-M} + \calF(\chi).
\end{align}
Since $u_n(t_n)$ tends to $\chi$ weakly in $\calX_2$ and
using once more Lemma~\ref{lemmalsc}, 
it is then easy to see that 
$\calF(u_n(t_n))$ tends to $\calF(\chi)$, which readily
entails that $u_n(t_n)\to\chi$ {\sl strongly}\/
in $\calX_2$, i.e.~{\sl (L3)}.\dimbox
\beos\label{ancheinfty}
 We point out that the attractor $\calA$ turns out to be
 more regular. More precisely, it is bounded
 and hence ``weakly'' compact in $\calXi$. Indeed, it is 
 easy to realize that the set of stationary points of (P)
 mentioned in property {\sl (L2)}\/ is also bounded in $\calXi$. 
 Moreover, \eqref{regoui2.2} entails that $\calS$ is (sequentially)
 ``weakly'' compact, i.e.~{\sl (L3)}\/ holds, in $\calXi$.
 Thus, Ball's procedure sketched in the Appendix can be
 repeated with respect to the ``weak'' topology 
 in $\calXi$. As a further consequence, 
 $\calA$ is also {\sl strongly}\/
 compact in $W^{2,p}(\Omega)$ for all $p\in[1,\infty)$.
\eddos
\beos\label{ancheasso}
 On account of the previous Remark, our procedure
 entails existence of an absorbing set $B_0$ for $\calS$
 which is bounded in $\calXi$ (not just in $\calX_2$).
\eddos


\section{Exponential attractors}
\label{secexpo}

In this section we prove Theorem \ref{teoattrattoreesp} 
by means of the method of $\ell$-trajectories. In
order to apply the theory of \cite{MP}
sketched in Section~\ref{preli}, we take 
$\calX:=V$ endowed with its standard norm. In comparison
with the global attractor, which was constructed in 
the smaller space $\calX_2$, we are thus working 
with weaker norm and topology (cf.~Remark~\ref{VnonX}
for additional comments on this point).

We know from the previous Section that $\calS$ 
admits an absorbing set $B_0$ bounded in $\calXi$.
We let (uniqueness holds on $B_0$,
thus we can use the ``semigroup'' $S(\cdot)$)
\beeq{defiB1}
   B_1:=\overline{\cup_{t\in[0,T_0]}S(t)B_0},
\end{equation}
where $T_0>0$ is such that $S(t)B_0\subset B_0$ 
for all $t\ge T_0$ and
the closure is taken w.r.t.~the {\sl weak}\/
topology of $\calXi$.
Due to the uniform character of estimate \eqref{regoui2.2}
(now the initial data are in $B_0$, so they are uniformly
bounded in $\calXi$), $B_1$ is still absorbing and 
bounded in $\calXi$. Moreover, we claim that
$B_1$ is positively invariant. To prove this fact, 
we let $\tau>0$ and assume that 
$u_0\in B_1$ is given by 
\beeq{coseu0}
  u_0=\lim_{n\nearrow\infty} S(t_n)u\zzn,
\end{equation}
where $\{u\zzn\}\subset B_0$ and $\{t_n\}\subset [0,T_0]$.
Then, using uniform boundedness,
weak compactness arguments and uniqueness
of solutions it is not difficult 
to realize that
\beeq{limiu0}
  S(t_n+\tau) u\zzn = S(\tau) \big(S(t_n)u\zzn\big)
   \to S(\tau)u_0
\end{equation}
{\sl weakly}\/ in $\calXi$ as $n\nearrow\infty$ 
(note that we cannot use 
directly~{\sl (S5)}\/ since we do not know that
$S(t_n)u\zzn$ converges {\sl strongly}\/ in $\calX_2$).
This readily entails that $S(\tau)u_0\in B_1$, which
is then positively invariant.

At this point, possibly making a 
positive and finite time shift,
we consider elements of $\calS$ starting
from initial data in $B_1$. Following \cite[Sec.~2]{MP}
and Section~\ref{preli} in this paper,
we set $\X_\ell:=L^{2}(0,\ell;\calX)$, where the choice
of $\ell\in(0,\infty)$ is here arbitrary, and define 
$\calBl^1$ as the set of $\ell$-trajectories 
whose initial datum lies in $B_1$. Using that
$B_1$ is positively invariant and {\sl weakly}\/
closed in $\calXi$, it is not difficult to 
show that $\calBl^1$ is a compact set in $\calX_\ell$,
so that, in particular, \eqref{calBl1} holds.
%

We now show the validity of conditions {\sl (M1)}, {\sl (M2)} and
{\sl (M3)}\/ reported in Section~\ref{preli}. To do this,
we prove a number of a priori estimates involving 
the difference of two solutions. Namely, we take $u_1,u_2$ 
solving (P) and starting from $u\zzu,u\zzd\in B_1$, respectively,
and set $u:=u_1-u_2$. Then, writing \eqref{eqn} for
$u=u_1$ and for $u=u_2$, and taking the difference, we have
\beeq{eqndiff}
 \alpha(\uut)-\alpha(\udt)+B u+W'(u_1)-W'(u_2)=0.
\end{equation}
In the sequel, the varying constant $c>0$ and the constants
$c_1,c_2,\dots>0$, whose numeration is restarted, will be allowed
to depend on $B_1$ and on $\ell$, additionally.
Let us test \eqref{eqndiff} by $u_t$.
We get
\beeq{contod11}
  \sigma\|u_t\|^2
   +\ddt\|u\|_V^2
   \le c\|u\|^2,
\end{equation}
where we also used the Young inequality and that,
thanks to \eqref{regoui2.2},
there exists $c>0$ depending on $B_1$ such that 
$\|W''(u_1(r))\|_\infty+\|W''(u_2(r))\|_\infty\le c$
for all $r\in[0,\infty)$. Then, by Gronwall's Lemma, 
\beeq{contod12}
  \|u(y)\|_V^2
   \le e^{c(y-s)}\|u(s)\|_V^2
   \le e^{2c\ell}\|u(s)\|_V^2
   =: c_1\|u(s)\|_V^2
\end{equation}
for all $s,y$ such that $0\le y-s\le 2\ell$.
Then, taking $s\in[0,\ell]$, $t\in[s,2\ell]$ and integrating
\eqref{contod11} over $[s,t]$, we infer
\beeq{contod13}
  \sigma\int_s^t\|u_t(r)\|^2\,\dir
   +\|u(t)\|_V^2
   \le c\int_s^t\|u(r)\|^2
   +\|u(s)\|_V^2.
\end{equation}
Thus, using \eqref{contod12} integrated for 
$y\in[s,t]$ to estimate the first term in the 
\rhs\ above, we get, for $t=2\ell$,
\beeq{contod14}
  \sigma\int_s^{2\ell}\|u_t(r)\|^2\,\dir
   +\|u(2\ell)\|_V^2
   \le c_2\|u(s)\|_V^2,
\end{equation}
whence, integrating for $s\in[0,\ell]$, 
\beeq{contod15}
  \sigma\ell\|u_t\|_{L^2(\ell,2\ell;H)}^2
   +\ell\|u(2\ell)\|_V^2
   \le c_2\|u\|_{L^2(0,\ell;V)}^2.
\end{equation}
Now, let us notice that a direct comparison argument
in the difference of the \eqref{eqn} written
for $u_1$ and for $u_2$ gives
\beeq{contod16}
  \|u\|_{H^{2}(\Omega)}^2
  \le c\big(\|u\|^2+\|Bu\|^2\big)
   \le c_3\|u\|^2 
   +c_3\|u_t\|^2,
\end{equation}
where the last inequality holds by the local Lipschitz 
continuity of $\alpha$ and $W'$ and the uniform
$\calXi$-boundedness of $u_1$, $u_2$.
Thus, evaluating the above formula in $y\in[\ell,2\ell]$, 
and using~\eqref{contod12},
\beeq{contod17}
  \|u(y)\|_{H^2(\Omega)}^2
   \le c_3c_1\|u(s)\|_V^2 
   +c_3\|u_t(y)\|^2.
\end{equation}
Finally, integrating for $s\in[0,\ell]$ and $y\in[\ell,2\ell]$
and recalling~\eqref{contod15},
\beeq{contod18}
  \|u\|_{L^2(\ell,2\ell;H^2(\Omega))}^2
   \le c_4\|u\|_{L^2(0,\ell;V)}^2.
\end{equation}
We are in the position to show properties {\sl (M1)}, {\sl (M2)} and
{\sl (M3)}. Setting
\beeq{defiWell}
  W_\ell:=\left\{v\in L^2(0,\ell;H^{2}(\Omega)):~v_t\in
    L^2(0,\ell;H)\right\},
\end{equation}
from \eqref{contod18} and \eqref{contod15} we have, respectively,
\begin{align}
  \|L_\ell u_1-L_\ell u_2\|_{L^{2}(0,\ell;H^{2}(\Omega))}
   \le c\|u_1-u_2\|_{L^{2}(0,\ell;V)},\label{contod19}\\
  \big\|(L_\ell u_1-L_\ell u_2)_t\big\|_{L^{2}(0,\ell;H)}
   \le c\|u_1-u_2\|_{L^{2}(0,\ell;V)}\label{contod20},
\end{align}
which imply property {\sl (M1)}\/ thanks to a straightforward 
application of the Aubin-Lions compactness Lemma.

Concerning {\sl (M2)}, this follows from \eqref{contod12} by
taking $y=s+t$, with $t$ varying in $[0,\tau]$, $\tau>0$, and 
integrating for $s\in[0,\ell]$ (the constant $c_1$
will actually take the value $e^{2c\tau}$, instead of 
$e^{2c\ell}$, with these choices).

Finally, property {\sl (M3)}\/ is a simple and
direct consequence of the time-regularity \eqref{stimaL2H}
of time derivatives of solutions 
(cf.~\cite[Lemma 2.2]{MP}).

According now to \cite[Theorem 2.5]{MP}, our procedure
entails existence of an exponential attractor $\mathcal{M}_\ell$
in the space of short trajectories. To show the 
existence of an exponential attractor also in the physical 
state space, we have to check the regularity
{\sl (M4)}\/ for the evaluation map $e$, which 
follows easily from \eqref{contod12}
by taking $y=\ell$ and integrating for $s\in[0,\ell]$.
Thus, thanks also to Remark~\ref{nuovaespo},
the set $\mathcal{M}:=e(\mathcal{M}_\ell)$
is an exponential attractor in $\calX=V$ for the semiflow~$\calS$.
\beos
 We stress once more that $\calM$ is a compact set in 
 $\calV$ (cf., however, Remark~\ref{VnonX}), but 
 it is able to attract exponentially fast 
 only the sets which are bounded in $\calX_2$ (and
 not all bounded sets if $V$).
\eddos




\section{Appendix}

We show here that the construction of global 
attractors for {\sl generalized}\/ semiflows 
(i.e., in our terminology, semiflows with ``strong-strong'' 
continuity properties but with no uniqueness at
all) given in \cite{Ba1} can be extended to our
situation. Actually, in comparison with J.~Ball's proof,
we have some simplification (mainly of 
technical character) due to the unique 
continuation~{\sl (S3)}. On the other hand, 
since our property~{\sl (S5)} is 
weaker than J.~Ball's ``strong-strong'' continuity
\cite[\sl (H4)]{Ba1}, we have to suitably modify
some points, which become now slightly more complicated.
For the reader's convenience we report
at least the highlights of all steps of J.~Ball's argument.
Concerning the proofs, we just point out the different 
points, instead. Basically, we will see that when in J.~Ball's proofs
\cite[\sl (H4)]{Ba1} occurs, we can replace it by
the combined use of~{\sl (S5)}\/ and the {\sl asymptotic
compactness}~{\sl (L3)}. In agreement with our
specific situation, the phase space will be indicated
as $\calX_2$ in what follows, but of course everything
holds for a generic metric space additionally endowed
with some ``weak'' topology.
\begin{prop}[Lemma 3.4 in \cite{Ba1}]\label{ball1}
\begin{sl}
 Let {\sl (S1)--(S5)}\/ and {\sl (L3)}\/ hold and let
 $B\subset \calX_2$ a bounded set. Then, the $\omega$-limit
 $\omega(B)$ is nonempty, compact, {\rm fully} invariant 
 and attracts $B$. 
\end{sl}
\end{prop}
\begin{proof}
 It is obvious from {\sl (L3)}\/
 that $\omega(B)$ is nonempty and easy to show
 directly that it is closed. We now prove that,
 for all $z\in\omega(B)$, there exists
 a {\sl complete trajectory}\/ $\psi$ taking values
 in $\omega(B)$ and such that $\psi(0)=z$
 (we recall that ``complete trajectory'' means
 that $\psi:\RR\to\calX_2$ is such that 
 $\psi(\cdot+\tau)\in \calS$ for all $\tau\in\RR$).
 Let then $\{u_n\}\subset\calS$ and $t_n\nearrow\infty$ such
 that $u_n(t_n)\to z$ and $\{u_n(0)\}\subset B$.
 By {\sl (S2)}, the sequence $\{v_n\}$, defined by
 $v_n(\cdot):=u_n(t_n+\cdot)$, lies in $\calS$
 and satisfies 
 $v_n(0)\to z$ strongly. Then,
 by {\sl (S5)}, there exist a nonrelabelled
 subsequence of $n$ and a solution 
 $v\in \calS$ such that, for all $t>0$,
 $u_n(t_n+t)=v_n(t)\to v(t)$ {\sl weakly}\/ in $\calX_2$.
 On the other hand, setting 
 $w_n(\cdot):=u_n(t_n/2+\cdot)$, it is $w_n\in\calS$.
 Moreover, we notice that, with no modifications in the proof, 
 it is still valid here \cite[Prop.~3.1]{Ba1}, which
 says that {\sl (L3)}\/ entails {\sl eventual boundedness},
 i.e., that for any bounded $B$ there exists $\tau_B\ge 0$
 such that $\cup_{t\ge\tau_B} T(t)B$ is still bounded
 (recall that $T(t)$ was defined in \eqref{defiTt}).
 Thus, we have that $\{w_n(0)\}$ is bounded and
 consequently, thanks to {\sl (L3)},
 $u_n(t_n+t)=w_n(t_n/2+t)$ converges
 strongly to its limit which is already identified
 as $v(t)$. Moreover,
 it is clear that $v(t)\in \omega(B)$ for all $t\ge 0$.
 This shows that from $z$ originates a (semi)trajectory $v$
 taking values in $\omega(B)$. The same trick used 
 above permits to adapt also J.~Ball's proof that
 $v$ extends to a complete trajectory $\psi$.
 Next, noting that on $\omega(B)$ uniqueness holds, 
 the above property also entails
 the {\sl complete invariance}\/
 of $\omega(B)$ (which did not necessarily
 hold in Ball's case). Finally, 
 the proof that $\omega(B)$ is compact and 
 attracts $B$ is essentially the
 same as in \cite{Ba1}. 
\end{proof}
\begin{prop}[Lemma 3.5 in \cite{Ba1}]\label{ball2}
\begin{sl}
 Let {\sl (S1)--(S5)}\/ and {\sl (L3)}\/ hold and let
 $\calS$ be {\rm pointwise dissipative}, namely let there 
 exist $B_0$ bounded in $\calX_2$ such that any $u\in\calS$
 eventually takes values in $B_0$. Then, there exist
 $\tau,\delta>0$ such that the set 
 \beeq{defiinto}
   B_1:=\bigcup_{t\ge\tau}T(t)(B(B_0,\delta)),
 \end{equation}
 with $B(B_0,\delta)$ denoting the open $\delta$-neighbourhood
 of $B_0$, is a bounded absorbing set for $\calS$.
\end{sl}
\end{prop}
\begin{proof}
 Let $\delta>0$. Then, by eventual boundedness, 
 there exists $\tau>0$ such that $B_1$ defined
 in \eqref{defiinto} is bounded.
 By contradiction, let us assume that some bounded $B$ 
 is not absorbed by $B_1$. Then, there exist
 $\{u_n\}\subset\calS$ and $t_n\nearrow\infty$ with
 $\{u_n(0)\}\subset B$ such that, for all $n$,
 $u_n(t_n)\not\in B_1$. Let us then
 set $v_n(\cdot):=u_n(t_n/2+\cdot)$, 
 so that $v_n(0)=u_n(t_n/2)$ and $v_n(t_n/2)=u_n(t_n)$.
 By {\sl (L3)}, at least for a subsequence,
 $v_n(0)\to z$ strongly. This entails by {\sl (S5)}\/
 that there exists $v\in \calS$ such that 
 $v_n(t)\to v(t)$ weakly for all $t\in[0,\infty)$.
 As before, since $v_n(t)=u_n(t_n/2+t)$ and 
 $\{u_n(0)\}$ is bounded, by {\sl (L3)}\/
 the convergence $v_n(t)\to v(t)$ is actually strong.
 Moreover, it is easy to see
 (proceed exactly as in \cite{Ba1})
 that $v_n(t)\not\in B(B_0,\delta)$
 for all $t\in[0,t_n/2-\tau]$. Thus, passing to 
 the (strong) limit, we have that 
 $v(t)\not\in B(B_0,\delta)$
 for all $t\in[0,\infty)$. Since $v$ 
 is a trajectory, this contradicts the point
 dissipativity of $\calS$ and gives the assert.
\end{proof}
\begin{prop}[Theorem 3.3 in \cite{Ba1}]\label{ball3}
\begin{sl}
 Let {\sl (S1)--(S5)}\/ and {\sl (L3)}\/ hold and let
 $\calS$ be {\rm pointwise dissipative}. Then, $\calS$ admits
 the global attractor $\calA$.
\end{sl}
\end{prop}
\begin{proof}
 It is as in \cite{Ba1}, up to minor modifications.
\end{proof}
\begin{prop}[Theorem 5.1 in \cite{Ba1}]\label{ball4}
\begin{sl}
 Let {\sl (S1)--(S5)}\/ and {\sl (L1)--(L3)}\/ hold.
 Then, $\calS$ is {\rm pointwise dissipative} (hence, 
 by the previous result, it admits the global attractor).
\end{sl}
\end{prop}
\begin{proof}
 Although it is similar to that in \cite{Ba1}, we
 prefer to give some more detail. First, it is easy to
 prove that, noting as $V$ the Liapounov functional
 and as $\calE_0$ the set of rest (i.e., stationary)
 points of $\calS$,
 given $u\in \calS$, $V$ is constant on $\omega(u)$ 
 and $\omega(u)$ is contained in $\calE_0$. To conclude, 
 we show that, given an arbitrary $\delta>0$, 
 any $u\in\calS$ eventually takes values in the 
 (bounded) set $B_0:=B(\calE_0,\delta)$. 
 Actually, if by contradiction $u(t_n)\not\in B_0$ for 
 a diverging sequence $\{t_n\}$,
 defining $v_n(\cdot):=u(t_n/2+\cdot)$
 and being, as before, $\{v_n\}\subset\calS$ and 
 $\{v_n(0)\}$ bounded, by asymptotic compactness
 $u(t_n)=v_n(t_n/2)$ has a subsequence which converges
 to an element of $\calE_0$.
\end{proof}
%




\vspace{10mm}

\large\noindent%
{\bf First author's address:}\\[1mm]
Giulio Schimperna\\
Dipartimento di Matematica, Universit\`a degli Studi di Pavia\\
Via Ferrata, 1,~~I-27100 Pavia,~~Italy\\
E-mail:~~{\tt giusch04@unipv.it}

\vspace{6mm}

\large\noindent%
{\bf Second author's address:}\\[1mm]
Antonio Segatti\\
Weierstrass Institute for Applied Analysis and Stochastics\\
Mohrenstrasse, 39,~~D-10117 Berlin,~~Germany\\
E-mail:~~{\tt segatti@wias-berlin.de}



\end{document}